\documentclass[10pt]{amsart}

\usepackage{hyperref}
\hypersetup{nesting=true,debug=true,naturalnames=true}
\usepackage{graphicx,amssymb,upref}
\usepackage[usenames,dvipsnames]{pstricks}
\usepackage{epsfig}
\usepackage{color}
\usepackage[utf 8]{inputenc}
\usepackage[T1]{fontenc}
\usepackage{amsmath, amsthm}
\usepackage[english]{babel}
\usepackage{amssymb}
\usepackage{latexsym}
\usepackage{graphicx,amssymb,upref,mathtools}
\usepackage[all]{xy}
\usepackage{tikz-cd}
\usepackage{yhmath}
\usetikzlibrary{matrix} 

\date{\today}

\newtheorem{theorem}{Theorem}[section]
\newtheorem{lemma}[theorem]{Lemma}

\newtheorem{corollary}[theorem]{Corollary}
\theoremstyle{definition}
\newtheorem{definition}[theorem]{Definition}
\newtheorem{example}[theorem]{Example}

\newtheorem{remark}[theorem]{Remark}

\newcommand{\ot}{\otimes}
\newcommand{\co}{\circ}

\newcommand{\IY}{I(q_Y)}
\newcommand{\IA}{I(q_A)}
\newcommand{\ID}{I(q_D)}
\newcommand{\IAp}{I(q_{A^{\prime}})}
\newcommand{\IDp}{I(q_{D^{\prime}})}

\DeclareMathOperator{\Hom}{Hom}

\title[Relative Rota-Baxter operators, modules and projections]{Relative Rota-Baxter operators, modules and projections}

\begin{document}
	
\maketitle

\begin{center}
{\bf Jos\'e Manuel Fern\'andez Vilaboa$^{1}$, Ram\'on
Gonz\'{a}lez Rodr\'{\i}guez$^{2}$ and Brais Ramos P\'erez$^{3}$}.
\end{center}
	
\vspace{0.4cm}
\begin{center}
{\small $^{1}$ [https://orcid.org/0000-0002-5995-7961].}
\end{center}
\begin{center}
{\small  CITMAga, 15782 Santiago de Compostela, Spain.}
\end{center}
\begin{center}
{\small  Universidade de Santiago de Compostela. Departamento de Matem\'aticas,  Facultade de Matem\'aticas, E-15771 Santiago de Compostela, Spain. 
\\ email: josemanuel.fernandez@usc.es.}
\end{center}
\vspace{0.2cm}
	
\begin{center}
{\small $^{2}$ [https://orcid.org/0000-0003-3061-6685].}
\end{center}
\begin{center}
{\small  CITMAga, 15782 Santiago de Compostela, Spain.}
\end{center}
\begin{center}
{\small  Universidade de Vigo, Departamento de Matem\'{a}tica Aplicada II,  E. E. Telecomunicaci\'on,
E-36310  Vigo, Spain.
\\email: rgon@dma.uvigo.es}
\end{center}
\vspace{0.2 cm}

\begin{center}
{\small $^{3}$ [https://orcid.org/0009-0006-3912-4483].}
\end{center}
\begin{center}
{\small  CITMAga, 15782 Santiago de Compostela, Spain. \\}
\end{center}
\begin{center}
{\small  Universidade de Santiago de Compostela. Departamento de Matem\'aticas,  Facultade de Matem\'aticas, E-15771 Santiago de Compostela, Spain. 
\\email: braisramos.perez@usc.es}
\end{center}
\vspace{0.2cm}
	
	
\begin{abstract}
The present article is devoted to introduce, in a braided monoidal setting, the notion of module over a relative Rota-Baxter operator. It is proved that there exists an adjunction between the category of modules associated to an invertible relative Rota-Baxter operator and the category of modules associated to a Hopf brace, which induces an equivalence by assuming certain additional hypothesis. Moreover, the notion of projection between relative Rota-Baxter operators is defined, and it is proved that those which are called ``strong'' give rise to a module according to the previous definition in the cocommutative setting.
\end{abstract} 

\vspace{0.2cm}



\section*{Introduction}
An important issue in the field of mathematical physics consists on finding solutions of the Quantum Yang-Baxter Equation (QYBE). The QYBE appeared in the 1970s in the field of quantum and statistical mechanics (see \cite{Yang} and \cite{Baxter}), and a complete classification of its solutions has not yet been obtained at this stage. A solution of such equation is an automorphism $c\colon V\otimes V\rightarrow V\otimes V$, where $V$ is a vector space over the field $\mathbb{F}$ and $\otimes$ denotes the tensor product of vector spaces over ${\mathbb F}$, which satisfies the following equality:
\begin{equation}\label{QYBE}\tag{QYBE}
(c\otimes id_{V})\circ(id_{V}\otimes c)\circ(c\otimes id_{V})=(id_{V}\otimes c)\circ(c\otimes id_{V})\circ (id_{V}\otimes c).
\end{equation}

Applications of the QYBE are diverse and cover different areas of mathematics and physics (to mention a few: knot theory, non-commutative geometry or quantum groups, between others), what implies that its solutions have been studied hardly and from different points of view. 

The first method to obtain solutions sistematically of the QYBE was proposed by Drinfeld in \cite{DR2}. He introduced the notion of quasitriangular Hopf algebra and proved that their associated modules give rise to solutions of the QYBE. It is an important result that, if $H$ is a Hopf algebra, the category of modules over its Drinfeld's double, $D(H)$, which is a quasitriangular Hopf algebra, is equivalent to the category of Yetter-Drinfeld modules over $H$ (see \cite[Theorem IX.5.2. and Section XIII.5]{K}). So, as a consequence of these facts, every Yetter-Drinfeld module over a Hopf algebra $H$ induce a solution of such equation. 

Later, Radford, Majid and Bespalov in \cite{RAD}, \cite{MAJ2} and \cite{Besp}, respectively, stated that, given a Hopf algebra $H$ such that its antipode is an isomorphism, there exists a categorical equivalence between  the category of Hopf algebras in ${}^{H}_{H}{\sf YD}$ and the category of Hopf algebra projections over $H$,
where ${}^{H}_{H}{\sf YD}$ denotes the category of Yetter-Drinfeld modules over $H$ and a Hopf algebra projection over $H$ is a pair of Hopf algebra morphisms, $
f\colon H\rightarrow B$ and $g\colon B\rightarrow H$, satisfying that $g\circ f=id_{H}.$ Therefore, taking into account the above mentioned results, to construct a solution of the QYBE it is enough to have a Hopf algebra projection.

On the other hand, in \cite{DR1} Drinfeld proposed to focus on the task of obtained set-theoretical solutions of the QYBE, which are those where a solution $c$ is a linear map induced by a mapping $\overline{c}\colon X\times X\rightarrow X\times X$, where $X$ is a set (in this situation, $V$ is the $\mathbb{F}$-vector space spanned by $X$). The study of this kind of solutions was pursued subsequently by several authors, for example, Etingof, Schedler and Soloviev in \cite{ESS} or Gateva-Ivanova in \cite{GI}. With this aim, the notion of brace was introduced by Rump in \cite{Rump} and then generalized by Guarnieri and Vendramin in \cite{GV} for the non-abelian setting obtaining the concept of skew brace. A  skew brace is a pair of groups $(G,\star)$ and $(G,\circ)$ which satisfies the following compatibility condition:
\begin{equation}\label{skewbrace}
g\circ(h\star t)=(g\circ h)\star g^{-1}\star(g\circ t),
\end{equation}
for all $g,h,t\in G$ and where $g^{-1}$ denotes the inverse of $g$ with regard to the group structure $(G,\star)$. In such paper the authors obtain that every skew brace induce a set-theoretical solution of the QYBE not always involutive, i.e., the inverse of such solution $c$ is not necessarily $c$. 

The quantum version of a skew brace is what is known by a Hopf brace, objects introduced by Angiono, Galindo and Vendramin in \cite{AGV}. Then, a Hopf brace is a pair of Hopf algebras sharing the underlying coalgebra structure, $H_{1}=(H, 1,\cdot, \epsilon,\Delta,\lambda)$ and $H_{2}=(H,1_{\circ},\circ,\epsilon,\Delta,S)$, which satisfy the following compatibility condition between the products:
\begin{equation}\label{hopfbrace}
g\circ(h\cdot t)=(g_{1}\circ h)\cdot\lambda(g_{2})\cdot(g_{3}\circ t),
\end{equation}
for all $g,h,t\in H$. As it occurs with skew braces, the subclass of the cocommutative Hopf braces also induce solutions of the QYBE (see \cite[Corollary 2.4]{AGV}).

Recently, by combining the two previously mentioned approaches to obtaining solutions for the QYBE, Fernández Vilaboa et al. in \cite{FGRR} studied the theory of projections in the Hopf brace setting, introducing the notion of Hopf brace projection, a suitable definition for the category of Yetter-Drinfeld modules over a Hopf brace and extending the correspondence of Radford-Majid-Bespalov to this framework.

Since the appearence of Hopf braces, another structures related with them have emerged. The first we will mention are invertible 1-cocycles, which appeared initially in Angiono's et al. article \cite{AGV}. Mixing \cite[Theorem 1.12]{AGV} and \cite[Theorem 3.2]{GONROD}, it was proved that the categories of invertible 1-cocycles and Hopf braces are equivalent. Also, Brzeziński in \cite{BRZ1} generalized the Hopf brace structure by modifying \eqref{hopfbrace} using a cocycle, what gave rise to the notion of Hopf truss. Besides, recently in \cite{LST}, Li, Sheng and Tang introduced the categories of post-Hopf algebras and relative Rota-Baxter operators. Regarding post-Hopf algebras, Fernández Vilaboa, González Rodríguez and Ramos Pérez proved in \cite{FGR} that the category of Hopf braces and a certain subcategory of post-Hopf algebras are isomorphic in the cocommutative setting. On the other hand, relative Rota-Baxter operators are a generalization of the notion of Rota-Baxter operator given by Goncharov in \cite{Goncharov} for cocommutative Hopf algebras. In \cite[Theorem 3.3]{LST} it was proved that there exists an adjunction between the category of post-Hopf algebras and the category of relative Rota-Baxter operators in a cocommutative context and, as a consequence of the above mentioned isomorphism between Hopf braces and post-Hopf algebras, this adjunction also holds between the category of Hopf braces and the category of relative Rota-Baxter operators under cocommutativity.

Looking at all the background, it is not unreasonable to tackle the study of projections for the objects mentioned in the previous paragraphs and whose connection with Hopf braces is strong. So, this paper is devoted to study the theory of modules and projections for relative Rota-Baxter operators.

The paper is organized as follows: After an initial preliminary section, where we are going to fix the notation and remember the basic necessary notions for the development of the article, Section \ref{sect2} is devoted to introduce the concept of module over a relative Rota-Baxter operator (Definition \ref{module_RB}). In this section, first we will prove that such category is symmetric monoidal under cocommutativity conditions and assuming that the base category is symmetric (Theorem \ref{symmonrRB}). Then, some functorial results are proved in order to state a correspondence between the category of modules over a certain Hopf brace and the category of modules over a relative Rota-Baxter operator. Using such correspondence, in Theorem \ref{adjVU} it is shown that there exists an adjunction between the category of modules over an invertible relative Rota-Baxter operator and the category of modules over such Hopf brace, which induce an equivalence by assuming some additional conditions (see Theorem \ref{equivVU}). Section \ref{sect3} is devoted to introduce the category of projections for relative Rota-Baxter operators (Definition \ref{defprojrB}). As happens between the category of Hopf braces and relative Rota-Baxter operators (see Theorem \ref{adjunction}), in Thereom \ref{projadjuncQR} we will see that there also exists an adjunction between their respective projection categories. To conclude the paper, we will define what a strong projection of relative Rota-Baxter operators is and we will prove that every such projection in the cocommutative setting gives rise to a module in the sense of Definition \ref{module_RB}, as happens in the classical theory of Hopf algebra projections.

\section{Preliminaries}\label{sect1}

From now on {\sf C} denotes a strict braided monoidal category with tensor product $\ot$, unit object $K$ and braiding $c$. Considering that it is well known that every non-strict monoidal category is monoidal equivalent to a strict one, we can assume without loss of generality that the category {\sf C} is strict and then we omit  explicitly  the associativity and unit constraints. Thus, the  results proved  in this paper for objects and morphisms in {\sf C} remain valid for every non-strict braided monoidal category which would include, for example, the category ${\mathbb F}$-{\sf Vect} of vector spaces over a field ${\mathbb F}$,  the category $R$-{\sf Mod} of left modules  over a commutative ring $R$, or the category {\sf Set} of sets. If for all $M,N\in {\sf C}$ the braiding satisfies that $c_{N,M}\circ c_{M,N}=id_{M\ot N}$, where $id$ denotes the identity morphism, we will say that ${\sf C}$ is symmetric. In what follows, for simplicity of notation, given objects $M$, $N$, $P$ in 
${\sf C}$ and a morphism $f:M\rightarrow N$, we write $P\ot f$ for
$id_{P}\ot f$ and $f \ot P$ for $f\ot id_{P}$.

\begin{definition}
An algebra in ${\sf  C}$ is a triple $A=(A, \eta_{A},
\mu_{A})$ where $A$ is an object in the category ${\sf  C}$ and
$\eta_{A}:K\rightarrow A$ (unit), $\mu_{A}:A\otimes A
\rightarrow A$ (product) are morphisms in ${\sf  C}$ such that
$\mu_{A}\circ (A\otimes \eta_{A})=id_{A}=\mu_{A}\circ
(\eta_{A}\otimes A)$, $\mu_{A}\circ (A\otimes \mu_{A})=\mu_{A}\circ
(\mu_{A}\otimes A)$. Given two algebras $A= (A, \eta_{A}, \mu_{A})$ and $B=(B, \eta_{B}, \mu_{B})$, a morphism  $f:A\rightarrow B$ in {\sf  C} is an algebra morphism if $\mu_{B}\circ (f\otimes f)=f\circ \mu_{A}$, $ f\circ
\eta_{A}= \eta_{B}$. 
		
If  $A$, $B$ are algebras in ${\sf  C}$, the tensor product
$A\otimes B$ is also an algebra in
${\sf  C}$ where
$\eta_{A\otimes B}\coloneqq\eta_{A}\otimes \eta_{B}$ and $\mu_{A\otimes B}\coloneqq(\mu_{A}\otimes \mu_{B})\circ (A\otimes c_{B,A}\otimes B).$
	
A coalgebra  in ${\sf  C}$ is a triple ${D} = (D,
\varepsilon_{D}, \delta_{D})$ where $D$ is an object in ${\sf
C}$ and $\varepsilon_{D}: D\rightarrow K$ (counit),
$\delta_{D}:D\rightarrow D\otimes D$ (coproduct) are morphisms in
${\sf  C}$ such that $(\varepsilon_{D}\otimes D)\circ
\delta_{D}= id_{D}=(D\otimes \varepsilon_{D})\circ \delta_{D}$,
$(\delta_{D}\otimes D)\circ \delta_{D}=
(D\otimes \delta_{D})\circ \delta_{D}.$ If ${D} = (D, \varepsilon_{D},
\delta_{D})$ and
${ E} = (E, \varepsilon_{E}, \delta_{E})$ are coalgebras, a morphism $f:D\rightarrow E$ in  {\sf  C} is a coalgebra morphism if $(f\otimes f)\circ
\delta_{D} =\delta_{E}\circ f$, $\varepsilon_{E}\circ f
=\varepsilon_{D}.$ 
		
Given  $D$, $E$ coalgebras in ${\sf  C}$, the tensor product $D\otimes E$ is a coalgebra in ${\sf  C}$ where $\varepsilon_{D\otimes E}\coloneqq\varepsilon_{D}\otimes \varepsilon_{E}$ and $\delta_{D\otimes E}\coloneqq(D\otimes c_{D,E}\otimes E)\circ( \delta_{D}\otimes \delta_{E}).$
\end{definition}

\begin{definition}
Let ${D} = (D, \varepsilon_{D},
\delta_{D})$ be a coalgebra and $A=(A, \eta_{A}, \mu_{A})$ an
algebra in $\sf{C}$. By $\operatorname{Hom}(D,A)$ we denote the set of morphisms
$f:D\rightarrow A$ in ${\sf  C}$. With the convolution operation
$f\ast g= \mu_{A}\circ (f\otimes g)\circ \delta_{D}$, $\operatorname{Hom}(D,A)$ is an algebra where the unit element is $\eta_{A}\circ \varepsilon_{D}=\varepsilon_{D}\otimes \eta_{A}$.
\end{definition}

\begin{definition}
Let  $A$ be an algebra. The pair
$(M,\varphi_{M})$ is a left $A$-module if $M$ is an object in
${\sf  C}$ and $\varphi_{M}:A\otimes M\rightarrow M$ is a morphism in ${\sf  C}$ satisfying $\varphi_{M}\circ(
\eta_{A}\ot M)=id_{M}$, $\varphi_{M}\circ (A\ot \varphi_{M})=\varphi_{M}\circ
(\mu_{A}\ot M)$. Given two left ${A}$-modules $(M,\varphi_{M})$
and $(N,\varphi_{N})$, $f:M\rightarrow N$ is a morphism of left ${A}$-modules if $\varphi_{N}\circ (A\ot f)=f\circ \varphi_{M}$. We will denote the category of left $A$-modules by ${}_{A}{\sf Mod}$.

Let  $D$ be a coalgebra. The pair
$(M,\rho_{M})$ is a left $D$-comodule if $M$ is an object in ${\sf  C}$ and $\rho_{M}:M\rightarrow D\ot M$ is a morphism
in ${\sf  C}$ satisfying $(\varepsilon_{D}\ot M)\co \rho_{M}=id_{M}$, $(D\ot \rho_{M})\co \rho_{M}=(\delta_{D}\ot M)\co \rho_{M}$. Given two left ${D}$-comodules $(M,\rho_{M})$
and $(N,\rho_{N})$, $f:M\rightarrow N$ is a morphism of left  ${D}$-comodules if $(D\ot f)\co \rho_{M}=\rho_{N}\co f$. We will denote the category of left $D$-comodules by $\;^{D}{\sf Comod}$.

In a similar way we can define the notions of right ${A}$-module and right $D$-comodule.  
\end{definition}
\begin{definition}
We say that $X$ is a
bialgebra  in ${\sf  C}$ if $(X, \eta_{X}, \mu_{X})$ is an
algebra, $(X, \varepsilon_{X}, \delta_{X})$ is a coalgebra, and
$\varepsilon_{X}$ and $\delta_{X}$ are algebra morphisms
(equivalently, $\eta_{X}$ and $\mu_{X}$ are coalgebra morphisms). Moreover, if there exists a morphism $\lambda_{X}:X\rightarrow X$ in ${\sf  C}$,
called the antipode of $X$, satisfying that $\lambda_{X}$ is the inverse of $id_{X}$ in $\operatorname{Hom}(X,X)$, i.e., 
\begin{equation}
\label{antipode}
id_{X}\ast \lambda_{X}= \eta_{X}\circ \varepsilon_{X}= \lambda_{X}\ast id_{X},
\end{equation}
we say that $X$ is a Hopf algebra. A morphism of Hopf algebras is an algebra-coalgebra morphism. Note that, if $f:X\rightarrow Y$ is a Hopf algebra morphism the following equality holds:
\begin{equation}
\label{morant}
\lambda_{Y}\co f=f\co \lambda_{X}.
\end{equation} 
		
With the composition of morphisms in {\sf C} we can define a category whose objects are Hopf algebras  and whose morphisms are morphisms of Hopf algebras. We denote this category by ${\sf  Hopf}$.
		
A Hopf algebra is commutative if $\mu_{X}\co c_{X,X}=\mu_{X}$ and cocommutative if $c_{X,X}\co \delta_{X}=\delta_{X}.$ In both cases $\lambda_{X}\circ \lambda_{X} =id_{X}$ and also, by \cite[Corollary 5]{Sch},  the identity 
\begin{equation}
\label{ccb}
c_{X,X}\circ c_{X,X}=id_{X\otimes X} 
\end{equation}
holds.
\end{definition}

If $X$ is a Hopf algebra, we have the following relevant properties of its antipode $\lambda_{X}$:  It is antimultiplicative and anticomultiplicative 
\begin{gather}
\label{a-antip1}
\lambda_{X}\co \mu_{X}=\mu_{X}\co (\lambda_{X}\ot \lambda_{X})\co c_{X,X},\\\label{a-antip2} \delta_{X}\co \lambda_{X}=c_{X,X}\co (\lambda_{X}\ot \lambda_{X})\co \delta_{X}, 
\end{gather}
and leaves the unit and counit invariant, i.e., 
\begin{gather}
\label{u-antip1}
\lambda_{X}\co \eta_{X}=  \eta_{X},\\\label{u-antip2} \varepsilon_{X}\co \lambda_{X}=\varepsilon_{X}.
\end{gather}
So, it is a direct consequence of these identities that, if $X$ is commutative, then $\lambda_{X}$ is an algebra morphism and, if $X$ is cocommutative, then $\lambda_{X}$ is a coalgebra morphism.

In the following definitions we recall the notion of left (co)module (co)algebra.

\begin{definition}
Let $X$ be a Hopf algebra. An algebra $A$  is said to be a left $X$-module algebra if $(A, \varphi_{A})$ is a left $X$-module and $\eta_{A}$, $\mu_{A}$ are morphisms of left $X$-modules, i.e.,
\begin{align}
\label{mod-alg1}
\varphi_{A}\circ (X\otimes \eta_{A})=\varepsilon_{X}\otimes \eta_{A},\\\label{mod-alg2}\varphi_{A}\circ (X\otimes \mu_{A})=\mu_{A}\circ \varphi_{A\otimes A},
\end{align}
where  $\varphi_{A\otimes A}\coloneqq(\varphi_{A}\otimes \varphi_{A})\circ (X\otimes c_{X,A}\otimes A)\circ (\delta_{X}\otimes A\otimes A)$ is the left action on $A\otimes A$. 
		
On the other hand,  $A$ is said to be a left $X$-comodule algebra if $(A,\rho_{A})$ is a left $X$-comodule and $\eta_{A}$ and $\mu_{A}$ are morphisms of left $X$-comodules, i.e.,
\begin{gather}
\label{comod-alg1} \rho_{A}\circ \eta_{A}=\eta_{X}\otimes \eta_{A},\\\label{comod-alg2} \rho_{A}\circ \mu_{A}=(X\otimes \mu_{A})\circ \rho_{A\otimes A}
\end{gather}
where $\rho_{A\otimes A}\coloneqq(\mu_{X}\ot A\ot A)\co (H\ot c_{A,X}\ot A)\co (\rho_{A}\ot \rho_{A})$ is the coaction on $A\ot A$. Equivalently, $(A,\rho_{A})$ is a left $X$-comodule algebra if and only if $\rho_{A}$ is an algebra morphism.
\end{definition}

\begin{definition}
Let $X$ be a Hopf algebra. A coalgebra $D$ is said to be a left $X$-module coalgebra if $(D,\varphi_{D})$ is a left $X$-module and $\varepsilon_{D}$, $\delta_{D}$ are morphisms of left $X$-modules, in other words, the following equalities hold:
\begin{gather}\label{mod-coalg1}
\varepsilon_{D}\circ\varphi_{D}=\varepsilon_{X}\otimes\varepsilon_{D},\\\label{mod-coalg2}\delta_{D}\circ\varphi_{D}=\varphi_{D\otimes D}\circ(X\otimes\delta_{D}).
\end{gather}
Equivalently, $(D,\varphi_{D})$ is a left $X$-module coalgebra if and only if $\varphi_{D}$ is a coalgebra morphism.
	
Finally, a coalgebra $D$  is said to be a left $X$-comodule coalgebra if $(D,\rho_{D})$ is a left $X$-comodule and  $\varepsilon_{D}$ and $\delta_{D}$ are morphisms of left $X$-comodules, i.e.,
\begin{gather}
\label{comod-coalg1}
(X\otimes \varepsilon_{D})\circ \rho_{D}=\eta_{X}\otimes \varepsilon_{D},\\\label{comod-coalg2}(X\otimes \delta_{D})\circ \rho_{D}=\rho_{D\otimes D}\circ \delta_{D}.
\end{gather}
\end{definition}
\begin{example}
Every Hopf algebra $X$ in {\sf C} has a structure of left module algebra over itself with the so called adjoint action $$\varphi_{X}^{ad}\coloneqq \mu_{X}\circ(\mu_{X}\otimes\lambda_{X})\circ(X\otimes c_{X,X})\circ(\delta_{X}\otimes X).$$
If $X$ is also cocommutative, then $(X,\varphi_{X}^{ad})$ is a left $X$-module algebra-coalgebra. Moreover, $X$ is a left $X$-comodule coalgebra with the adjoint coaction $\rho_{X}^{ad}\coloneqq (\mu_{X}\ot X)\co (X\ot c_{X,X})\co (\delta_{X}\ot \lambda_{X})\co \delta_{X}$.
\end{example}

\section{Relative Rota-Baxter operators and their modules}\label{sect2}
The present section of this paper is devoted to introduce what a module over a relative Rota-Baxter operator is. Relative Rota-Baxter operators have been introduced considering the underlying category ${\sf C}=\mathbb{F}\text{-}{\sf Vect}$ by Li, Sheng and Tang in \cite{LST} as a generalization of Rota-Baxter operators defined by Goncharov in \cite{Goncharov} for cocommutative Hopf algebras.

First of all we start by remembering the notion and basic properties of relative Rota-Baxter operators, as well as the strong relationship between these structures and Hopf braces, and we show that the category formed by this objects is symmetric monoidal. 

After that we will focus on the study of modules over a relative Rota-Baxter operator, giving a definition that allow us to show that there exists an adjunction between the category of modules over a Hopf brace (using the definition of these objects introduced by González in \cite{RGON}) and the category of modules over an invertible relative Rota-Baxter operator assuming cocommutativity. The section will finish by proving that, under certain additional hypothesis, the previous adjunction gives rise to an equivalence of categories.
\begin{definition}\label{rRB}
Let $H=(H,\eta_{H},\mu_{H},\varepsilon_{H},\delta_{H},\lambda_{H})$ and $B=(B,\eta_{B},\mu_{B},\varepsilon_{B},\delta_{B},\lambda_{B})$ be Hopf algebras in $\sf{C}$ such that $(H,\varphi_{H})$ is a left $B$-module algebra-coalgebra. We will say that a coalgebra morphism $T\colon H\rightarrow B$ is a relative Rota-Baxter operator if the following condition holds:
\begin{equation}\label{rRBcond}
\mu_{B}\circ (T\otimes T)=T\circ \mu_{H}\circ (H\otimes(\varphi_{H}\circ(T\otimes H)))\circ(\delta_{H}\otimes H).
\end{equation}
In what follows we will denote relative Rota-Baxter operators by $\left(T\begin{array}{c}H\\\downarrow\\B\end{array},\varphi_{H}\right)$.
	
If $\left(T\begin{array}{c}H\\\downarrow\\B\end{array},\varphi_{H}\right)$ and $\left(L\begin{array}{c}A\\\downarrow\\D\end{array},\varphi_{A}\right)$ are relative Rota-Baxter operators, a morphism between them is a pair $(f,h)$ where $f\colon H\rightarrow A$ and $h\colon B\rightarrow D$ are Hopf algebra morphisms and the following conditions hold:
\begin{gather}\label{cond1morrRB}
L\circ f=h\circ T,\\
\label{cond2morrRB}
f\circ\varphi_{H}=\varphi_{A}\circ(h\otimes f).
\end{gather}
	
Considering the natural composition of morphisms, relative Rota-Baxter operators give rise to a category that we will denote by {\sf rRB}. Moreover, we will denote by ${\sf rRB}^{\star}$ to the full subcategory of ${\sf rRB}$ whose objects are relative Rota-Baxter operators $\left(T\begin{array}{c}H\\\downarrow\\B\end{array},\varphi_{H}\right)$ such that $H$ is cocommutative, and by ${\sf coc}\textnormal{-}{\sf rRB}$ to the full subcategory of ${\sf rRB}^{\star}$ satisfying that both Hopf algebras involved, $H$ and $B$, are cocommutative.
\end{definition}

An important property of relative Rota-Baxter operators is that they preserve the unit. This will be proven in the following result.

\begin{lemma}
If $\left(T\begin{array}{c}H\\\downarrow\\B\end{array},\varphi_{H}\right)$ is a relative Rota-Baxter operator, then \begin{equation}
\label{etaT}
\eta_{B}=T\circ\eta_{H}.
\end{equation}
\end{lemma}
\begin{proof}
By \eqref{rRBcond}, the condition of morphism of left $B$-modules for $\eta_{H}$, the condition of coalgebra morphism for $T$ and the (co)unit property, we obtain that
\begin{align}\label{proof1etaT}
\mu_{B}\circ((T\circ\eta_{H})\otimes(T\circ\eta_{H}))=T\circ\eta_{H}.
\end{align}
Then, we have that
\begin{align*}
&\eta_{B}\\=&\eta_{B}\circ\varepsilon_{B}\circ T\circ\eta_{H}\;\footnotesize\textnormal{(by the condition of coalgebra morphism for $T$ and (co)unit properties)}\\=&\mu_{B}\circ(\lambda_{B}\otimes B)\circ\delta_{B}\circ T\circ \eta_{H}\;\footnotesize\textnormal{(by \eqref{antipode} for $B$)}\\=&\mu_{B}\circ(\lambda_{B}\otimes B)\circ(T\otimes T)\circ(\eta_{H}\otimes\eta_{H})\;\footnotesize\textnormal{(by the condition of coalgebra morphism for $T$ and $\eta_{H}$)}\\=&\mu_{B}\circ(\lambda_{B}\otimes B)\circ((T\circ\eta_{H})\otimes(\mu_{B}\circ((T\circ\eta_{H})\otimes(T\circ\eta_{H}))))\;\footnotesize\textnormal{(by \eqref{proof1etaT})}\\=&\mu_{B}\circ((\mu_{B}\circ(\lambda_{B}\otimes B)\circ\delta_{B}\circ T\circ\eta_{H})\otimes(T\circ\eta_{H}))\;\footnotesize\textnormal{(by associativity of $\mu_{B}$ and the condition of}\\&\footnotesize\textnormal{coalgebra morphism $T$ and $\eta_{H}$)}\\=&T\circ\eta_{H}\;\footnotesize\textnormal{(by \eqref{antipode} for $B$, the condition of coalgebra morphism for $T$ and (co)unit properties)}.\qedhere
\end{align*}
\end{proof}

If ${\sf C}$ is symmetric, ${\sf rRB}$ admits an structure of  symmetric monoidal category as can be seen in what follows.

\begin{theorem}\label{monoidalrRB}
Let's assume {\sf C} to be symmetric. The category of relative Rota-Baxter operators is a strict symmetric monoidal with tensor functor
\begin{align*}
\otimes\colon {\sf rRB}\times{\sf rRB}&\longrightarrow {\sf rRB}\\\left(\left(T\begin{array}{c}H\\\downarrow\\B\end{array},\varphi_{H}\right),\left(L\begin{array}{c}A\\\downarrow\\D\end{array},\varphi_{A}\right)\right)&\longmapsto \left(T\otimes L\begin{array}{c}H\otimes A\\\downarrow\\B\otimes D\end{array},\varphi_{H\otimes A}^t \right),
\end{align*}
where $\varphi_{H\otimes A}^t\coloneqq(\varphi_{H}\otimes\varphi_{A})\circ(B\otimes c_{D,H}\otimes A)$, unit object $\left(id_{K}\begin{array}{c}K\\\downarrow\\K\end{array},id_{K}\right)$ and symmetry given by
\begin{align*}
\tau_{T,L}\coloneqq (c_{H,A},c_{B,D})\colon \left(T\otimes L\begin{array}{c}H\otimes A\\\downarrow\\B\otimes D\end{array},\varphi_{H\otimes A}^t\right)\rightarrow \left(L\otimes T\begin{array}{c}A\otimes H\\\downarrow\\D\otimes B\end{array},\varphi_{A\otimes H}^t\right).
\end{align*}
\end{theorem}
\begin{proof}
When {\sf C} is symmetric, if $H=(H,\eta_{H},\mu_{H},\varepsilon_{H},\delta_{H},\lambda_{H})$ and $A=(A,\eta_{A},\mu_{A},\varepsilon_{A},\delta_{A},\lambda_{A})$ are Hopf algebras in {\sf C}, then $H\otimes A=(H\otimes A,\eta_{H}\otimes\eta_{A},\mu_{H\otimes A},\varepsilon_{H}\otimes\varepsilon_{A},\delta_{H\otimes A},\lambda_{H}\otimes\lambda_{A})$ is also a Hopf algebra in {\sf C}. 
	
Moreover, $(H\otimes A,\varphi_{H\otimes A})$ is a left $B\otimes D$-module algebra-coalgebra. Indeed, on the one side, the left module axioms are straightforward thanks to naturality of $c$. On the other side, it results easy to prove that $\varphi_{H\otimes A}^t\circ(B\otimes D\otimes \eta_{H}\otimes\eta_{A})=\varepsilon_{B}\otimes\varepsilon_{D}\otimes\eta_{H}\otimes\eta_{A}$ and the condition of morphism of left $B\otimes D$-modules for $\mu_{H\otimes A}$ follows by
\begin{align*}
&\mu_{H\otimes A}\circ(\varphi_{H\otimes A}^t\otimes\varphi_{H\otimes A}^t)\circ(B\otimes D\otimes c_{B\otimes D,H\otimes A}\otimes H\otimes A)\circ(\delta_{B\otimes D}\otimes H \otimes A\otimes H\otimes A)\\=&((\mu_{H}\circ(\varphi_{H}\otimes\varphi_{H})\circ(B\otimes c_{B,H}\otimes H)\circ(\delta_{B}\otimes H\otimes H))\otimes(\mu_{A}\circ(\varphi_{A}\otimes\varphi_{A})\circ(D\otimes c_{D,A}\otimes A)))\\&\circ(B\otimes ((H\otimes c_{D,H}\otimes D\otimes A)\circ(c_{D,H}\otimes c_{D,H}\otimes A)\circ(D\otimes c_{D,H}\otimes H\otimes A)\circ(\delta_{D}\otimes H\otimes c_{A,H}))\otimes A)\\&\footnotesize\textnormal{(by naturality of $c$ and {\sf C} symmetric)}\\=&((\mu_{H}\circ(\varphi_{H}\otimes\varphi_{H})\circ(B\otimes c_{B,H}\otimes H)\circ(\delta_{B}\otimes H\otimes H))\otimes(\mu_{A}\circ(\varphi_{A}\otimes \varphi_{A})\circ(D\otimes c_{D,A}\otimes A)\\&\circ(\delta_{D}\otimes A\otimes A)))\circ(B\otimes((H\otimes c_{D,H}\otimes A)\circ(c_{D,H}\otimes c_{A,H}))\otimes A)\;\footnotesize\textnormal{(by naturality of $c$)}\\=&((\varphi_{H}\circ(B\otimes\mu_{H}))\otimes(\varphi_{A}\circ(D\otimes\mu_{A})))\circ(B\otimes((H\otimes c_{D,H}\otimes A)\circ(c_{D,H}\otimes c_{A,H}))\otimes A)\;\footnotesize\textnormal{(by the condition of}\\&\footnotesize\textnormal{morphism of $B$-modules for $\mu_{H}$ and the condition of morphism of $D$-modules for $\mu_{A}$)}\\=&\varphi_{H\otimes A}^t\circ(B\otimes D\otimes\mu_{H\otimes A})\;\footnotesize\textnormal{(by naturality of $c$)}. 
\end{align*}
In addition, $(\varepsilon_{H}\otimes\varepsilon_{A})\circ\varphi_{H\otimes A}^t=\varepsilon_{B}\otimes\varepsilon_{D}\otimes\varepsilon_{H}\otimes\varepsilon_{A}$ is straightforward while $$\delta_{H\otimes A}\circ\varphi_{H\otimes A}^t=(\varphi_{H\otimes A}^t\otimes\varphi_{H\otimes A}^t)\circ(B\otimes D\otimes c_{B\otimes D,H\otimes A}\otimes H\otimes A)\circ(\delta_{B\otimes D}\otimes\delta_{H\otimes A})$$ follows by
\begin{align*}
&\delta_{H\otimes A}\circ\varphi_{H\otimes A}^{t}\\=&(H\otimes c_{H,A}\otimes A)\circ(((\varphi_{H}\otimes\varphi_{H})\circ(B\otimes c_{B,H}\otimes H)\circ(\delta_{B}\otimes\delta_{H}))\otimes((\varphi_{A}\otimes\varphi_{A})\circ(D\otimes c_{D,A}\otimes A)\\&\circ(\delta_{D}\otimes\delta_{A})))\circ(B\otimes c_{D,H}\otimes A)\;\footnotesize\textnormal{(by the condition of coalgebra morphism for $\varphi_{H}$ and $\varphi_{A}$)}\\=&(\varphi_{H\otimes A}^t\otimes\varphi_{H\otimes A}^{t})\circ(B\otimes D\otimes c_{B\otimes D,H\otimes A}\otimes H\otimes A)\circ(\delta_{B\otimes D}\otimes\delta_{H\otimes A})\;\footnotesize\textnormal{(by naturality of $c$ and {\sf C} symmetric)}.
\end{align*}
	
Also note that, due to be considering the standard structure of tensor coproduct, if $T$ and $L$ are coalgebra morphism, then $T\otimes L$ is also a coalgebra morphism. Therefore, to conclude the monoidal character of {\sf rRB} we only have to compute that \eqref{rRBcond} holds. Indeed,
\begin{align*}
&(T\otimes L)\circ\mu_{H\otimes A}\circ(H\otimes A\otimes(\varphi_{H\otimes A}^t\circ(T\otimes L\otimes H\otimes A)))\circ(\delta_{H\otimes A}\otimes H\otimes A)\\=&((T\circ\mu_{H}\circ (H\otimes(\varphi_{H}\circ(T\otimes H)))\circ(\delta_{H}\otimes H))\otimes (L\circ\mu_{A}\circ(A\otimes(\varphi_{A}\circ(L\otimes A)))\circ(\delta_{A}\otimes A)))\\&\circ(H\otimes c_{A,H}\otimes A)\;\footnotesize\textnormal{(by naturality of $c$ and {\sf C} symmetric)}\\=&((\mu_{B}\circ(T\otimes T))\otimes(\mu_{D}\circ(L\otimes L)))\circ(H\otimes c_{A,H}\otimes A)\;\footnotesize\textnormal{(by \eqref{rRBcond} for $T$ and $L$)}\\=&\mu_{B\otimes D}\circ((T\otimes L)\otimes(T\otimes L))\;\footnotesize\textnormal{(by naturality of $c$)}.
\end{align*}
	
On the other hand, $\tau_{T,L}$ is a symmetry for ${\sf rRB}$ because, when the base category ${\sf C}$ is symmetric, $(c_{H,A},c_{B,D})$ is a morphism in {\sf rRB} between $\left(T\otimes L\begin{array}{c}H\otimes A\\\downarrow\\B\otimes D\end{array},\varphi_{H\otimes A}^t\right)$ and $\left(L\otimes T\begin{array}{c}A\otimes H\\\downarrow\\D\otimes B\end{array},\varphi_{A\otimes H}^t\right)$.
\end{proof}

After proving these general properties which relative Rota-Baxter operators satisfy, we will see that under suitable conditions there exists a functorial link between  relative Rota-Baxter operators and Hopf braces. First we will remember the definition of Hopf brace and its main properties.

\begin{definition}
\label{H-brace}
Let $H=(H, \varepsilon_{H}, \delta_{H})$ be a coalgebra in {\sf C}. Let's assume that there are two algebra structures $(H, \eta_{H}^1, \mu_{H}^1)$, $(H, \eta_{H}^2, \mu_{H}^2)$ defined on $H$ and suppose that there exist two endomorphism of $H$ denoted by $\lambda_{H}^{1}$ and $\lambda_{H}^{2}$. We will say that 
$$(H, \eta_{H}^{1}, \mu_{H}^{1}, \eta_{H}^{2}, \mu_{H}^{2}, \varepsilon_{H}, \delta_{H}, \lambda_{H}^{1}, \lambda_{H}^{2})$$
is a Hopf brace in {\sf C} if:
\begin{itemize}
\item[(i)] $H_{1}=(H, \eta_{H}^{1}, \mu_{H}^{1},  \varepsilon_{H}, \delta_{H}, \lambda_{H}^{1})$ is a Hopf algebra in {\sf C}.
\item[(ii)] $H_{2}=(H, \eta_{H}^{2}, \mu_{H}^{2},  \varepsilon_{H}, \delta_{H}, \lambda_{H}^{2})$ is a Hopf algebra in {\sf C}.
\item[(iii)] The  following equality holds:
$$\mu_{H}^{2}\co (H\ot \mu_{H}^{1})=\mu_{H}^{1}\co (\mu_{H}^{2}\ot \Gamma_{H_{1}} )\co (H\ot c_{H,H}\ot H)\co (\delta_{H}\ot H\ot H),$$
\end{itemize}
where  $$\Gamma_{H_{1}}\coloneqq\mu_{H}^{1}\co (\lambda_{H}^{1}\ot \mu_{H}^{2})\co (\delta_{H}\ot H).$$
	
For any Hopf brace
$\eta_{H}^{1}=\eta_{H}^2$ holds and therefore, because of this property, the expression of a Hopf brace  is reduced to
$$(H, \eta_{H}, \mu_{H}^{1}, \mu_{H}^{2}, \varepsilon_{H}, \delta_{H}, \lambda_{H}^{1}, \lambda_{H}^{2}).$$
	
In the following lines, a Hopf brace will be denoted by ${\mathbb H}=(H_{1}, H_{2})$ or in a simpler way by $\mathbb{H}$.
	
\end{definition}

\begin{definition}
If  ${\mathbb H}$ is a Hopf brace in {\sf C}, we will say that ${\mathbb H}$ is cocommutative if $\delta_{H}=c_{H,H}\circ \delta_{H}$, i.e., if $H_{1}$ and $H_{2}$ are cocommutative Hopf algebras in {\sf C}.
\end{definition}

\begin{definition}
Given two Hopf braces ${\mathbb H}$  and  ${\mathbb D}$ in {\sf C}, a morphism $x$ in {\sf C} between the two underlying objects is called a morphism of Hopf braces if both $x:H_{1}\rightarrow D_{1}$ and $x:H_{2}\rightarrow D_{2}$ are Hopf algebra morphisms.
	
Hopf braces together with morphisms of Hopf braces form a category which we denote by {\sf HBr}. Moreover, cocommutative Hopf braces constitute a full subcategory of $\sf{HBr}$ which we will denote by ${\sf coc}\textnormal{-}{\sf HBr}$.
\end{definition}

Moreover, in our braided context \cite[Lemma 1.8]{AGV} and \cite[Remark 1.9]{AGV} hold and then we have that the algebra $(H,\eta_{H}, \mu_{H}^{1})$ is a left $H_{2}$-module algebra with action $\Gamma_{H_{1}}$ and $\mu_{H}^2$ admits the following expression:
\begin{equation}
\label{eb2}
\mu_{H}^2=\mu_{H}^{1}\circ (H\otimes \Gamma_{H_{1}})\circ (\delta_{H}\otimes H). 
\end{equation}
In addition, by \cite[Lemma 2.2]{AGV}, $\Gamma_{H_{1}}$ is a coalgebra morphism when $\mathbb{H}$ is cocommutative.

The following result is the braided monoidal version of the result  proved by Li et al. in \cite[Proposition 3.2 and Theorem 3.3]{LST} for Hopf braces in the category of vector spaces over a field $\mathbb{F}$. We do the proof in detail to clarify certain properties and notations that will be very useful in the rest of the paper.

\begin{theorem}\label{adjunction}
There exists a functor ${\sf F}\colon {\sf coc}\textnormal{-}{\sf HBr}\longrightarrow{\sf rRB}^{\star}$ defined on objects by $${\sf F}(\mathbb{H})=\left(id_{H}\begin{array}{c}H_{1}\\\downarrow\\H_{2}\end{array},\Gamma_{H_{1}}\right),$$  and on morphisms by ${\sf F}(x)=(x,x)$. 
	
Moreover, there exists a functor ${\sf G}\colon {\sf rRB}^{\star}\longrightarrow {\sf coc}\textnormal{-}{\sf HBr}$ defined on objects by $${\sf G}\left(\left(T\begin{array}{c}H\\\downarrow\\B\end{array},\varphi_{H}\right)\right)=\overline{\mathbb{H}},$$ where $\overline{\mathbb{H}}=(H,\overline{H})$ is the Hopf brace with  $\overline{H}=(H,\eta_{H},\overline{\mu}_{H},\varepsilon_{H},\delta_{H},\overline{\lambda}_{H})$  the Hopf algebra whose product and antipode are given by
\begin{gather*}
\overline{\mu}_{H}\coloneqq \mu_{H}\circ(H\otimes(\varphi_{H}\circ(T\otimes H)))\circ(\delta_{H}\otimes H),\\\overline{\lambda}_{H}\coloneqq\varphi_{H}\circ((\lambda_{B}\circ T)\otimes \lambda_{H})\circ\delta_{H},
\end{gather*}
and on morphisms by ${\sf G}(f,h)=f$.
	
In addition, $F$ factors by ${\sf coc}\textnormal{-}{\sf rRB}$ and ${\sf F}$ is left adjoint to ${\sf G}$.
\end{theorem}
\begin{proof}
First of all, let's see that ${\sf F}$ is well-defined. If $\mathbb{H}=(H_{1},H_{2})$ is a cocommutative Hopf brace in {\sf C}, then we already know that $(H_{1},\Gamma_{H_{1}})$ is a left $H_{2}$-module algebra and, thanks to the cocommutativity, $\Gamma_{H_{1}}$ is a coalgebra morphism, i.e., $(H_{1},\Gamma_{H_{1}})$ is a left $H_{2}$-module algebra-coalgebra. Moreover, \eqref{eb2} implies that \eqref{rRBcond} holds. Therefore, $\left(id_{H}\begin{array}{c}H_{1}\\\downarrow\\H_{2}\end{array},\Gamma_{H_{1}}\right)$ is a relative Rota-Baxter operator. In addition, if $\mathbb{D}$ is another cocommutative Hopf brace and $x\colon \mathbb{H}\rightarrow \mathbb{D}$ is a morphism of Hopf braces, then the pair $(x,x)$ is a morphism of relative Rota-Baxter operators between $\left(id_{H}\begin{array}{c}H_{1}\\\downarrow\\H_{2}\end{array},\Gamma_{H_{1}}\right)$ and $\left(id_{D}\begin{array}{c}D_{1}\\\downarrow\\D_{2}\end{array},\Gamma_{D_{1}}\right)$. Indeed, it is straightforward to compute that \eqref{cond1morrRB} holds and \eqref{cond2morrRB} follows by
\begin{align*}
&x\circ \Gamma_{H_{1}}\\\nonumber=&\mu_{D}^{1}\circ ((x\circ\lambda_{H}^{1})\otimes(x\circ\mu_{H}^{2}))\circ(\delta_{H}\otimes H)\;\footnotesize\textnormal{(by the condition of algebra morphism for $x\colon H_{1}\rightarrow D_{1}$)}\\\nonumber=&\mu_{D}^{1}\circ(\lambda_{D}^{1}\otimes\mu_{D}^{2})\circ(((x\otimes x)\circ\delta_{H})\otimes x)\;\footnotesize\textnormal{(by \eqref{morant} and the condition of algebra morphism for $x\colon H_{2}\rightarrow D_{2}$)}\\\nonumber=&\Gamma_{D_{1}}\circ(x\otimes x)\;\footnotesize\textnormal{(by the condition of coalgebra morphism for $x$)}.
\end{align*}
	
Now, let's proof that ${\sf G}$ is well-defined. On the one hand, consider $\left(T\begin{array}{c}H\\\downarrow\\B\end{array},\varphi_{H}\right)$ a relative Rota-Baxter operator such that $H$ is cocommutative. Let's show that $\overline{H}$ is a Hopf algebra. At first, note that it is straightforward to prove that $\eta_{H}$ is the unit for $\overline{\mu}_{H}$ and the associativity of $\overline{\mu}_{H}$ follows by
\begin{align*}
&\overline{\mu}_{H}\circ(\overline{\mu}_{H}\otimes H)\\=&\mu_{H}\circ (H\otimes (\varphi_{H}\circ(T\otimes H)))\circ (((\mu_{H}\otimes\mu_{H})\circ (H\otimes c_{H,H}\otimes H)\circ (\delta_{H}\otimes((\varphi_{H}\otimes\varphi_{H})\circ(B\otimes c_{B,H}\otimes H)\\&\circ(((T\otimes T)\circ\delta_{H})\otimes\delta_{H}))))\otimes H)\circ(\delta_{H}\otimes H\otimes H)\;\footnotesize\textnormal{(by the condition of coalgebra morphism for $\mu_{H}$, $\varphi_{H}$ and $T$)}\\=&\mu_{H}\circ(H\otimes\varphi_{H})\circ(\overline{\mu}_{H}\otimes(T\circ\overline{\mu}_{H})\otimes H)\circ(((H\otimes c_{H,H}\otimes H)\circ(\delta_{H}\otimes\delta_{H}))\otimes H)\;\footnotesize\textnormal{(by the cocommutativity}\\&\footnotesize\textnormal{and coassociativity of $\delta_{H}$ and naturality of $c$)}\\=&\mu_{H}\circ (\overline{\mu}_{H}\otimes(\varphi_{H}\circ((\mu_{B}\circ(T\otimes T))\otimes H)))\circ (((H\otimes c_{H,H}\otimes H)\circ(\delta_{H}\otimes\delta_{H}))\otimes H)\;\footnotesize\textnormal{(by \eqref{rRBcond})}\\=&\mu_{H}\circ(H\otimes(\mu_{H}\circ (\varphi_{H}\otimes\varphi_{H})\circ(B\otimes c_{B,H}\otimes H)\circ(\delta_{B}\otimes H\otimes H)))\circ(((H\otimes T)\circ\delta_{H})\\&\otimes((H\otimes (\varphi_{H}\circ(T\otimes H)))\circ (\delta_{H}\otimes H)))\;\footnotesize\textnormal{(by module axioms for $(H,\varphi_{H})$, coassociativity of $\delta_{H}$, the condition}\\&\footnotesize\textnormal{of coalgebra morphism for $T$ and associativity of $\mu_{H}$)}\\=&\overline{\mu}_{H}\circ(H\otimes\overline{\mu}_{H})\;\footnotesize\textnormal{(by the condition of morphism of left $B$-modules for $\mu_{H}$)}.
\end{align*}
Moreover, note that 
\begin{align*}
&\delta_{H}\circ\overline{\mu}_{H}\\=&(\mu_{H}\otimes\mu_{H})\circ(H\otimes c_{H,H}\otimes H)\circ(\delta_{H}\otimes((\varphi_{H}\otimes\varphi_{H})\circ(B\otimes c_{B,H}\otimes H)\circ (((T\otimes T)\circ\delta_{H})\otimes\delta_{H})))\circ(\delta_{H}\otimes H)\\&\footnotesize\textnormal{(by the condition of coalgebra morphism for $\mu_{H}$, $\varphi_{H}$ and $T$)}\\=&((\mu_{H}\circ(H\otimes(\varphi_{H}\circ(T\otimes H))))\otimes(\mu_{H}\circ(H\otimes(\varphi_{H}\circ(T\otimes H)))))\circ(H\otimes ((H\otimes c_{H,H}\otimes H)\\&\circ(c_{H,H}\otimes c_{H,H}))\otimes H)\circ(((\delta_{H}\otimes\delta_{H})\circ\delta_{H})\otimes\delta_{H})\;\footnotesize\textnormal{(by naturality of $c$)}\\=&(\overline{\mu}_{H}\otimes\overline{\mu}_{H})\circ(H\otimes c_{H,H}\otimes H)\circ(\delta_{H}\otimes\delta_{H})\;\footnotesize\textnormal{(by cocommutativity and coassociativity of $\delta_{H}$)}.
\end{align*}
and, by the condition of coalgebra morphism for $\mu_{H}$, $\varphi_{H}$ and $T$ and the counit property, the equality $\varepsilon_{H}\circ\overline{\mu}_{H}=\varepsilon_{H}\otimes\varepsilon_{H}$ also holds. Thus, to conclude that $\overline{H}$ is a Hopf algebra, it only remains to prove that $\overline{\lambda}_{H}$ is the inverse of $id_{H}$ for the convolution in $\Hom(H,\overline{H})$, operation that we will denote by $\overline{\ast}$. Firstly,
\begin{align*}\label{antipode1side}
&id_{H}\,\overline{\ast}\, \overline{\lambda}_{H}\\\nonumber=&\mu_{H}\circ (H\otimes(\varphi_{H}\circ (((id_{B}\ast \lambda_{B})\circ T)\otimes\lambda_{H})\circ\delta_{H}))\circ\delta_{H}\;\footnotesize\textnormal{(by the coassociativity of $\delta_{H}$, module axioms for $(H,\varphi_{H})$}\\\nonumber&\footnotesize\textnormal{and the condition of coalgebra morphism for $T$)}\\\nonumber=&id_{H}\ast\lambda_{H}\;\footnotesize\textnormal{(by \eqref{antipode} for $B$, the condition of coalgebra morphism for $T$, the counit property and module axioms for $(H,\varphi_{H})$)}\\\nonumber=&\varepsilon_{H}\otimes\eta_{H}\;\footnotesize\textnormal{(by \eqref{antipode} for $H$)}.
\end{align*}
Note also that $\overline{\lambda}_{H}$ satisfies the following property: $\overline{\lambda}_{H}$ is a coalgebra morphism because
\begin{align*}
&\delta_{H}\circ\overline{\lambda}_{H}\\\nonumber=&(\varphi_{H}\otimes\varphi_{H})\circ(B\otimes c_{B,H}\otimes H)\circ((\delta_{B}\circ\lambda_{B}\circ T)\otimes(\delta_{H}\circ\lambda_{H}))\circ\delta_{H}\;\footnotesize\textnormal{(by the condition of coalgebra morphism for $\varphi_{H}$)}\\\nonumber=&(\varphi_{H}\otimes\varphi_{H})\circ(B\otimes c_{B,H}\otimes H)\circ(((\lambda_{B}\otimes\lambda_{B})\circ(T\otimes T)\circ c_{H,H}\circ\delta_{H})\otimes((\lambda_{H}\otimes\lambda_{H})\circ c_{H,H}\circ\delta_{H}))\circ\delta_{H}\;\footnotesize\textnormal{(by \eqref{a-antip2},}\\\nonumber&\footnotesize\textnormal{the condition of coalgebra morphism for $T$ and naturality of $c$)}\\\nonumber=&((\varphi_{H}\circ((\lambda_{B}\circ T)\otimes\lambda_{H}))\otimes(\varphi_{H}\circ((\lambda_{B}\circ T)\otimes\lambda_{H})))\circ(H\otimes c_{H,H}\otimes H)\circ(\delta_{H}\otimes\delta_{H})\circ\delta_{H}\;\footnotesize\textnormal{(by cocommutativity}\\&\nonumber\footnotesize\textnormal{of $\delta_{H}$ and naturality of $c$)}\\\nonumber=&(\overline{\lambda}_{H}\otimes\overline{\lambda}_{H})\circ\delta_{H}\;\footnotesize\textnormal{(by coassociativity and cocommutativity of $\delta_{H}$)}
\end{align*}
and 
\begin{equation*}
\varepsilon_{H}\circ\overline{\lambda}_{H}=\varepsilon_{H}
\end{equation*}
what follows by the condition of coalgebra morphism for $\varphi_{H}$ and $T$, \eqref{u-antip2} and the counit property. As a consequence,
\begin{align*}
&\overline{\lambda}_{H}\circ\overline{\lambda}_{H}\\\nonumber=&\overline{\mu}_{H}\circ((\eta_{H}\circ\varepsilon_{H})\otimes(\overline{\lambda}_{H}\circ\overline{\lambda}_{H}))\circ\delta_{H}\;\footnotesize\textnormal{(by the (co)unit property)}\\\nonumber=&\overline{\mu}_{H}\circ((\overline{\mu}_{H}\circ(H\otimes\overline{\lambda}_{H})\circ\delta_{H})\otimes(\overline{\lambda}_{H}\circ\overline{\lambda}_{H}))\circ\delta_{H}\;\footnotesize\textnormal{(by $id_{H}\,\overline{\ast}\, \overline{\lambda}_{H}=\eta_{H}\otimes \varepsilon_{H}$)}\\\nonumber=&\overline{\mu}_{H}\circ(H\otimes (\overline{\mu}_{H}\circ (\overline{\lambda}_{H}\otimes(\overline{\lambda_{H}}\circ\overline{\lambda}_{H}))\circ\delta_{H}))\circ\delta_{H}\;\footnotesize\textnormal{(by coassociativity of $\delta_{H}$ and associativity of $\overline{\mu}_{H}$)}\\\nonumber=&\overline{\mu}_{H}\circ(H\otimes((id_{H}\,\overline{\ast}\,\overline{\lambda}_{H})\circ\overline{\lambda}_{H}))\circ\delta_{H}\;\footnotesize\textnormal{(by the condition of coalgebra morphism of $\overline{\lambda}_{H}$)}\\\nonumber=&id_{H}\;\footnotesize\textnormal{(by $id_{H}\,\overline{\ast}\, \overline{\lambda}_{H}=\varepsilon_{H}\otimes \eta_{H}$, the condition of coalgebra morphism of $\overline{\lambda}_{H}$ and the (co)unit property)}.
\end{align*}
Therefore,
\begin{align*}
&\overline{\lambda}_{H}\,\overline{\ast}\, id_{H}\\\nonumber=&\overline{\mu}_{H}\circ(\overline{\lambda}_{H}\otimes(\overline{\lambda}_{H}\circ\overline{\lambda}_{H}))\circ\delta_{H}\;\footnotesize\textnormal{(by $\overline{\lambda}_{H}\circ\overline{\lambda}_{H}=id_{H}$)}\\\nonumber=&(id_{H}\,\overline{\ast}\,\overline{\lambda}_{H})\circ\overline{\lambda}_{H}\;\footnotesize\textnormal{(by the condition of coalgebra morphism of $\overline{\lambda}_{H}$)}\\\nonumber=&\varepsilon_{H}\otimes\eta_{H}\;\footnotesize\textnormal{(by $id_{H}\,\overline{\ast}\, \overline{\lambda}_{H}=\varepsilon_{H}\otimes \eta_{H}$ and the condition of coalgebra morphism of $\overline{\lambda}_{H}$)}.
\end{align*}
So, to conclude that $\overline{\mathbb{H}}$ is a Hopf brace it is enough to see that (iii) of Definition \ref{H-brace} holds. Indeed, 
\begin{equation}
\label{overlGammaH}
\overline{\Gamma}_{H}=\varphi_{H}\circ(T\otimes H)
\end{equation}
what follows by
\begin{align*}
\label{overlGammaH}
&\overline{\Gamma}_{H}\\=&\mu_{H}\circ((\lambda_{H}\ast id_{H})\otimes(\varphi_{H}\circ(T\otimes H)))\circ(\delta_{H}\otimes H)\;\footnotesize\textnormal{(by the coassociativity of $\delta_{H}$ and associativity of $\mu_{H}$)}\\=&\varphi_{H}\circ(T\otimes H)\;\footnotesize\textnormal{(by \eqref{antipode} and the (co)unit property)}.
\end{align*}
Then, 
\begin{align*}
&\mu_{H}\circ(\overline{\mu}_{H}\otimes\overline{\Gamma}_{H})\circ(H\otimes c_{H,H}\otimes H)\circ(\delta_{H}\otimes H\otimes H)\\=&\mu_{H}\circ(H\otimes (\mu_{H}\circ (\varphi_{H}\otimes\varphi_{H})\circ(B\otimes c_{B,H}\otimes H)\circ((\delta_{B}\circ T)\otimes H\otimes H)))\circ(\delta_{H}\otimes H\otimes H)\;\footnotesize\textnormal{(by \eqref{overlGammaH},}\\&\footnotesize\textnormal{coassociativity of $\delta_{H}$, associativity of $\mu_{H}$, naturality of $c$ and the condition of coalgebra morphism for $T$)}\\=&\overline{\mu}_{H}\circ(H\otimes\mu_{H})\;\footnotesize\textnormal{(by the condition of morphism of left $B$-modules for $\mu_{H}$)}.
\end{align*}
On the other hand, if $(f,h)$ is a morphism in ${\sf rRB}^{\star}$ between $\left(T\begin{array}{c}H\\\downarrow\\B\end{array},\varphi_{H}\right)$ and $\left(L\begin{array}{c}A\\\downarrow\\D\end{array},\varphi_{A}\right)$, then $f$ is a morphism of Hopf braces between $\overline{\mathbb{H}}$ and $\overline{\mathbb{A}}$. Indeed,
\begin{align*}
&f\circ\overline{\mu}_{H}\\=&\mu_{A}\circ(f\otimes(f\circ\varphi_{H}\circ(T\otimes H)))\circ(\delta_{H}\otimes H)\;\footnotesize\textnormal{(by the condition of algebra morphism for $f$)}\\=&\mu_{A}\circ(f\otimes(\varphi_{A}\circ((h\circ T)\otimes f)))\circ(\delta_{H}\otimes H)\;\footnotesize\textnormal{(by \eqref{cond2morrRB})}\\=&\mu_{A}\circ(A\otimes (\varphi_{A}\circ (L\otimes A)))\circ (((f\otimes f)\circ\delta_{H})\otimes f)\;\footnotesize\textnormal{(by \eqref{cond1morrRB})}\\=&\overline{\mu}_{A}\circ(f\otimes f)\;\footnotesize\textnormal{(by the condition of coalgebra morphism for $f$)}.
\end{align*}

To prove that ${\sf F}$ is left adjoint to ${\sf G}$, it is enough to consider the bijection
\[{}^{\mathbb{H}}\Theta_{L}\colon\Hom_{{\sf coc}\textnormal{-}{\sf HBr}}(\mathbb{H}=(H_{1},H_{2}),\overline{\mathbb{A}}=(A,\overline{A}))\longrightarrow \Hom_{{\sf rRB}^{\star}}\left(\left(id_{H}\begin{array}{c}H_{1}\\\downarrow\\H_{2}\end{array},\Gamma_{H_{1}}\right),\left(L\begin{array}{c}A\\\downarrow\\D \end{array},\varphi_{A}\right) \right)\]
given by ${}^{\mathbb{H}}\Theta_{L}(y)=(y,L\circ y)$ and $({}^{\mathbb{H}}\Theta_{L})^{-1}(f,h)=f$ for all $\mathbb{H}\in {\sf coc}\textnormal{-}{\sf HBr}$ and $\left(L\begin{array}{c}A\\\downarrow\\D \end{array},\varphi_{A}\right)\in{\sf rRB}^{\star}$.
\end{proof}

Taking into account the previous considerations, in the following definition the notion of module over a relative Rota-Baxter operator is introduced.
\begin{definition}\label{module_RB}
	Let $\left(T\begin{array}{c}H\\\downarrow\\B\end{array},\varphi_{H}\right)$ be a relative Rota-Baxter operator. We will say that a 6-tuple $$(M,N,\phi_{M},\varphi_{M},\varphi_{N},\gamma)$$ is a left module over  $\left(T\begin{array}{c}H\\\downarrow\\B\end{array},\varphi_{H}\right)$ if the conditions:
	\begin{enumerate}
		\item[(i)] $(M,\phi_{M})$ is a left $H$-module and $(M,\varphi_{M})$ is a left $B$-module such that the equality
		\begin{equation}\label{compatmodHmodB}
			\varphi_{M}\circ(B\otimes\phi_{M})=\phi_{M}\circ(\varphi_{H}\otimes\varphi_{M})\circ(B\otimes c_{B,H}\otimes M)\circ(\delta_{B}\otimes H\otimes M)
		\end{equation}
		holds,
		\item[(ii)]  $(N,\varphi_{N})$ is a left $B$-module,
		\item[(iii)]  $\gamma\colon M\rightarrow N$ is a morphism satisfying that 
		\begin{equation}\label{compatmodT} \varphi_{N}\circ(T\otimes\gamma)=\gamma\circ\phi_{M}\circ (H\otimes(\varphi_{M}\circ(T\otimes M)))\circ(\delta_{H}\otimes M).\end{equation}
	\end{enumerate}
	
	Let $(M,N,\phi_{M},\varphi_{M},\varphi_{N},\gamma)$ and $(P,Q,\phi_{P},\varphi_{P},\varphi_{Q},\theta)$ be left modules over $\left(T\begin{array}{c}H\\\downarrow\\B\end{array},\varphi_{H}\right)$. We will say that a pair $(r,s)$ is a morphism of left modules over $\left(T\begin{array}{c}H\\\downarrow\\B\end{array},\varphi_{H}\right)$ between $(M,N,\phi_{M},\varphi_{M},\varphi_{N},\gamma)$ and $(P,Q,\phi_{P},\varphi_{P},\varphi_{Q},\theta)$ if the conditions
	\begin{enumerate}
		\item[(i)] $r\colon (M,\phi_{M})\rightarrow (P,\phi_{P})$ is a morphism of left $H$-modules,
		\item[(ii)] $r\colon (M,\varphi_{M})\rightarrow (P,\varphi_{P})$ is a morphism of left $B$-modules,
		\item[(iii)] $s\colon (N,\varphi_{N})\rightarrow (Q,\varphi_{Q})$ is a morphism of left $B$-modules,
		\item[(iv)] $s\circ\gamma=\theta\circ r$
	\end{enumerate} 
	hold. 
	
	Therefore, with the obvious composition of morphisms, left modules over the relative Rota-Baxter operator $\left(T\begin{array}{c}H\\\downarrow\\B\end{array},\varphi_{H}\right)$ constitute a category that we will denote by ${}_{(T,\varphi_{H})}{\sf Mod}$.
\end{definition}
\begin{example}
	Note that, given $\left(T\begin{array}{c}H\\\downarrow\\B\end{array},\varphi_{H}\right)$ a relative Rota-Baxter operator, then $(H,B,\mu_{H},\varphi_{H},\mu_{B},T)$ is an object in ${}_{(T,\varphi_{H})}{\sf Mod}$. Moreover, $(K,K,\varepsilon_{H},\varepsilon_{B},\varepsilon_{B},id_{K})$ is called the trivial left module over any relative Rota-Baxter operator $\left(T\begin{array}{c}H\\\downarrow\\B\end{array},\varphi_{H}\right)$.
\end{example}
\begin{theorem}\label{symmonrRB}
Let's assume that {\sf C} is symmetric and let $\left(T\begin{array}{c}H\\\downarrow\\B\end{array},\varphi_{H}\right)$ an object in ${\sf coc}\textnormal{-}{\sf rRB}$. Then, the category ${}_{(T,\varphi_{H})}{\sf Mod}$ is monoidal with unit object $(K,K,\varepsilon_{H},\varepsilon_{B},\varepsilon_{B},id_{K})$ and tensor functor defined by 
\begin{align*}
\otimes\colon{}_{(T,\varphi_{H})}{\sf Mod}\times {}_{(T,\varphi_{H})}{\sf Mod}&\longrightarrow {}_{(T,\varphi_{H})}{\sf Mod}\\((M,N,\phi_{M},\varphi_{M},\varphi_{N},\gamma),(P,Q,\phi_{P},\varphi_{P},\varphi_{Q},\theta))&\longmapsto (M\otimes P,N\otimes Q,\phi_{M\otimes P},\varphi_{M\otimes P},\varphi_{N\otimes Q},\gamma\otimes\theta).
\end{align*}
	
Moreover, ${}_{(T,\varphi_{H})}{\sf Mod}$ is symmetric with symmetry isomorphism given by 
\[\tau_{(M,N,\gamma), (P,Q,\theta)}\coloneqq (c_{M,P},c_{N,Q}).\]
\end{theorem}
\begin{proof}
Consider $(M,N,\phi_{M},\varphi_{M},\varphi_{N},\gamma)$ and $(P,Q,\phi_{P},\varphi_{P},\varphi_{Q},\theta)$ objects in ${}_{(T,\varphi_{H})}{\sf Mod}$. Then, lelt's see that $(M\otimes P,N\otimes Q,\phi_{M\otimes P},\varphi_{M\otimes P},\varphi_{N\otimes Q},\gamma\otimes\theta)$ is also an object in ${}_{(T,\varphi_{H})}{\sf Mod}$. Indeed, note that, due to the monoidal character of the module categories ${}_{H}{\sf Mod}$ and ${}_{B}{\sf Mod}$, the only facts that remain us to compute is that \eqref{compatmodHmodB} and \eqref{compatmodT} hold. At first, we have that
\begin{align*}
&\phi_{M\otimes P}\circ(\varphi_{H}\otimes \varphi_{M\otimes P})\circ(B\otimes c_{B,H}\otimes M\otimes P)\circ(\delta_{B}\otimes H\otimes M\otimes P)\\=&(\phi_{M}\otimes\phi_{P})\circ (H\otimes c_{H,M}\otimes P)\circ(((\varphi_{H}\otimes \varphi_{H})\circ (B\otimes c_{B,H}\otimes H)\circ(\delta_{B}\otimes\delta_{H}))\otimes((\varphi_{M}\otimes\varphi_{P})\\&\circ (B\otimes c_{B,M}\otimes P)\circ (\delta_{B}\otimes M\otimes P)))\circ (B\otimes c_{B,H}\otimes M\otimes P)\circ (\delta_{B}\otimes H\otimes M\otimes P)\;\footnotesize\textnormal{(by the condition of}\\&\footnotesize\textnormal{coalgebra morphism for $\varphi_{H}$)}\\=&((\phi_{M}\circ(\varphi_{H}\otimes\varphi_{M})\circ(B\otimes c_{B,H}\otimes M)\circ (\delta_{B}\otimes H\otimes M))\otimes(\phi_{P}\circ(\varphi_{H}\otimes\varphi_{P})\circ (B\otimes c_{B,H}\otimes P)\circ (\delta_{B}\otimes H\otimes P)))\\&\circ (B\otimes((H\otimes c_{B,M}\otimes H)\circ (c_{B,H}\otimes c_{H,M}))\otimes P)\circ (\delta_{B}\otimes\delta_{H}\otimes M\otimes P)\;\footnotesize\textnormal{(by naturality of $c$, symmetrical character}\\&\footnotesize\textnormal{of {\sf C} and cocommutativity of $\delta_{B}$)}\\=&((\varphi_{M}\circ (B\otimes\phi_{M}))\otimes(\varphi_{P}\circ(B\otimes\phi_{P})))\circ (B\otimes((H\otimes c_{B,M}\otimes H)\circ (c_{B,H}\otimes c_{H,M}))\otimes P)\circ (\delta_{B}\otimes\delta_{H}\otimes M\otimes P)\\&\footnotesize\textnormal{(by \eqref{compatmodHmodB})}\\=&\varphi_{M\otimes P}\circ(B\otimes\phi_{M\otimes P})\;\footnotesize\textnormal{(by naturality of $c$)},
\end{align*}
so \eqref{compatmodHmodB} holds and, on the other side,
\begin{align*}
&(\gamma\otimes\theta)\circ\phi_{M\otimes P}\circ(H\otimes(\varphi_{M\otimes P}\circ(T\otimes M\otimes P)))\circ (\delta_{H}\otimes M\otimes P)\\=&((\gamma\circ\phi_{M}\circ (H\otimes (\varphi_{M}\circ(T\otimes M))))\otimes(\theta\circ \phi_{P}\circ (H\otimes (\varphi_{P}\circ(T\otimes P)))))\circ(H\otimes((H\otimes c_{H,M}\otimes H)\\&\circ(c_{H,H}\otimes c_{H,M}))\otimes P)\circ(((\delta_{H}\otimes\delta_{H})\circ\delta_{H})\otimes M\otimes P)\;\footnotesize\textnormal{(by naturality of $c$ and the condition of coalgebra}\\&\footnotesize\textnormal{morphism for $T$)}\\=&((\gamma\circ\phi_{M}\circ (H\otimes(\varphi_{M}\circ(T\otimes M)))\circ(\delta_{H}\otimes M))\otimes(\theta\circ\phi_{P}\circ (H\otimes (\varphi_{P}\circ(T\otimes P)))\circ(\delta_{H}\otimes P)))\\&\circ(H\otimes c_{H,M}\otimes P)\circ(\delta_{H}\otimes M\otimes P)\;\footnotesize\textnormal{(by coassociativity and cocommutativity of $\delta_{H}$ and naturality of $c$)}\\=&((\varphi_{N}\circ (T\otimes\gamma))\otimes(\varphi_{Q}\circ (T\otimes\theta)))\circ (H\otimes c_{H,M}\otimes P)\circ (\delta_{H}\otimes M\otimes P)\;\footnotesize\textnormal{(by \eqref{compatmodT})}\\=&\varphi_{N\otimes Q}\circ (T\otimes \gamma\otimes\theta)\;\footnotesize\textnormal{(by naturality of $c$ and the condition of coalgebra morphism for $T$)},
\end{align*}
what implies that \eqref{compatmodT} holds.
	
On the other hand, the symmetric character of ${}_{(T,\varphi_{H})}{\sf Mod}$ follows by the fact that, when ${\sf C}$ is symmetric and $H$ and $B$ are cocommutative Hopf algebras, the pair $(c_{M,P},c_{N,Q})$ is a morphism in ${}_{(T,\varphi_{H})}{\sf Mod}$.
\end{proof}

\begin{theorem}	If $(f,h)\colon \left(T\begin{array}{c}H\\\downarrow\\B\end{array},\varphi_{H}\right)\rightarrow \left(L\begin{array}{c}A\\\downarrow\\D\end{array},\varphi_{A}\right)$ is a morphism of relative Rota-Baxter operators, then there exists a functor 
\[{\sf R}_{(f,h)}\colon {}_{(L,\varphi_{D})}{\sf Mod}\longrightarrow {}_{(T,\varphi_{H})}{\sf Mod}\]
acting on objects by 
\[{\sf R}_{(f,h)}((M,N,\phi_{M},\varphi_{M},\varphi_{N},\gamma))=(M,N,\phi_{M}^{T}\coloneqq \phi_{M}\circ(f\otimes M),\varphi_{M}^{T}\coloneqq \varphi_{M}\circ(h\otimes M),\varphi_{N}^{T}\coloneqq\varphi_{N}\circ(h\otimes N)),\gamma)\]
and on morphisms by the identity.  
\end{theorem}
\begin{proof}
Let's show that it is well-defined on objects. First of all note that it is straightforward to prove that $(M,\phi_{M}^{T})$ is a left $H$-module and $(M,\varphi_{M}^{T})$ and $(N,\varphi_{N}^{T})$ are left $B$-modules. Thus, it only remains to see that equalities \eqref{compatmodHmodB} and \eqref{compatmodT} hold. Then,
\begin{align*}
&\phi_{M}^{T}\circ (\varphi_{H}\otimes\varphi_{M}^{T})\circ (B\otimes c_{B,H}\otimes M)\circ(\delta_{B}\otimes H\otimes M)\\=&\phi_{M}\circ (\varphi_{A}\otimes\varphi_{M})\circ(D\otimes c_{D,A}\otimes M)\circ(\delta_{D}\otimes A\otimes M)\circ(h\otimes f\otimes M)\;\footnotesize\textnormal{(by \eqref{cond2morrRB} for $(f,h)$, naturality of $c$}\\&\footnotesize\textnormal{and the condition of coalgebra morphism for $h$)}\\=&\varphi_{M}^{T}\circ(B\otimes\phi_{M}^{T})\;\footnotesize\textnormal{(by \eqref{compatmodHmodB})},
\end{align*}
i.e., \eqref{compatmodHmodB} holds, and also
\begin{align*}
&\gamma\circ \phi_{M}^{T}\circ(H\otimes(\varphi_{M}^{T}\circ (T\otimes M)))\circ (\delta_{H}\otimes M)\\=&\gamma\circ \phi_{M}\circ (A\otimes (\varphi_{M}\circ(L\otimes M)))\circ (((f\otimes f)\circ\delta_{H})\otimes M)\;\footnotesize\textnormal{(by \eqref{cond1morrRB} for $(f,h)$)}\\=&\gamma\circ \phi_{M}\circ (A\otimes (\varphi_{M}\circ (L\otimes M)))\circ ((\delta_{A}\circ f)\otimes M)\;\footnotesize\textnormal{(by the condition of coalgebra morphism for $f$)}\\=&\varphi_{N}\circ ((L\circ f)\otimes\gamma)\;\footnotesize\textnormal{(by \eqref{compatmodT})}\\=&\varphi_{N}^{T}\circ (T\otimes\gamma)\;\footnotesize\textnormal{(by \eqref{cond1morrRB} for $(f,h)$)},
\end{align*} 
what implies that \eqref{compatmodT} holds.
\end{proof}

Let's recall the notion of module over a Hopf brace introduced in \cite{RGON} and some results related with this concept which can be consulted in \cite{FGRR}.
\begin{definition}
Let $\mathbb{H}=(H_{1},H_{2})$ be a Hopf brace. We will say that a triple $(M,\psi_{M}^{1},\psi_{M}^{2})$ is a left module over $\mathbb{H}$ if $(M,\psi_{M}^{1})$ is a left $H_{1}$-module, $(M,\psi_{M}^{2})$ is a left $H_{2}$-module and the following compatibility condition holds:
\begin{equation}\label{compatmodbrace}
\psi_{M}^{2}\circ(H\otimes\psi_{M}^{1})=\psi_{M}^{1}\circ(\mu_{H}^{2}\otimes\Gamma_{M})\circ(H\otimes c_{H,H}\otimes M)\circ(\delta_{H}\otimes H\otimes M),
\end{equation}
where $$\Gamma_{M}\coloneqq \psi_{M}^{1}\circ(\lambda_{H}^{1}\otimes\psi_{M}^{2})\circ(\delta_{H}\otimes M).$$
	
If $(M,\psi_{M}^{1},\psi_{M}^{2})$ and $(N,\psi_{N}^{1},\psi_{N}^{2})$ are modules over the Hopf brace $\mathbb{H}$ and $f\colon M\rightarrow N$ is a morphism between them, we will say that $f$ is a morphism of left $\mathbb{H}$-modules if $f$ is a morphism of left $H_{1}$-modules and left $H_{2}$-modules.
	
Then, with the obvious composition of morphisms, modules over $\mathbb{H}$ constitute a category that we will denote by ${}_{\mathbb{H}}{\sf Mod}$. 
\end{definition}

\begin{remark}
Let $(M,\psi_{M}^{1},\psi_{M}^{2})$ be a left module over a Hopf brace $\mathbb{H}$. Composing on the right hand side of \eqref{compatmodbrace} with $H\otimes\eta_{H}\otimes M$ we obtain that the equality
\begin{equation}\label{psi2expression}
\psi_{M}^{2}=\psi_{M}^{1}\circ(H\otimes \Gamma_{M})\circ(\delta_{H}\otimes M)
\end{equation}
holds.
\end{remark}

The proof of the following result can be seen in \cite[Lemma 2.11]{FGRR}.
\begin{theorem}
Let $\mathbb{H}$ be a Hopf brace and $(M,\psi_{M}^{1},\psi_{M}^{2})$ a left module over $\mathbb{H}$. The equality 
\begin{equation}\label{GammaMPsiM1}
\Gamma_{M}\circ (H\otimes\psi_{M}^{1})=\psi_{M}^{1}\circ (\Gamma_{H_{1}}\otimes\Gamma_{M})\circ (H\otimes c_{H,H}\otimes M)\circ (\delta_{H}\otimes H\otimes M)
\end{equation}
holds and $(M,\Gamma_{M})$ is a left $H_{2}$-module.
\end{theorem}

Taking into account the previous results, in what follows we are going to construct two functors that set a relationship between the category of modules over a Hopf brace and the category of modules over a relative Rota-Baxter operator.

\begin{theorem}
Let $\mathbb{H}$ be a cocommutative Hopf brace. There exists a functor 
\[{\sf W}_{\mathbb{H}}\colon {}_{\mathbb{H}}{\sf Mod}\longrightarrow {}_{(id_{H}, \Gamma_{H_{1}})}{\sf Mod},\]
where $(id_{H}, \Gamma_{H_{1}})$ denotes the relative Rota-Baxter operator ${\sf F}(\mathbb{H})=\left(id_{H}\begin{array}{c}H_1\\\downarrow\\H_2\end{array},\Gamma_{H_{1}}\right)$ introduced in Theorem \ref{adjunction}, which acts on objects by $${\sf W}_{\mathbb{H}}((M,\psi_{M}^{1},\psi_{M}^{2}))=(M,M,\psi_{M}^{1},\Gamma_{M},\psi_{M}^{2},id_{M})$$ and on morphisms by ${\sf W}_{\mathbb{H}}(f)=(f,f)$. 
\end{theorem}
\begin{proof}
At first we are going to prove that if $(M,\psi_{M}^{1},\psi_{M}^{2})$ is a left module the Hopf brace $\mathbb{H}$, then the 6-tuple $(M,M,\psi_{M}^{1},\Gamma_{M},\psi_{M}^{2},id_{M})$ is a left module over the relative Rota-Baxter operator ${\sf F}(\mathbb{H})$. Indeed, on the one hand, \eqref{compatmodHmodB} follows by \eqref{GammaMPsiM1} and, on the other hand, \eqref{compatmodT} follows by \eqref{psi2expression}. Therefore, ${\sf W}$ is well-defined on objects.
	
On the other hand, if $f\colon (M,\psi_{M}^{1},\psi_{M}^{2})\rightarrow (N,\psi_{N}^{1},\psi_{N}^{2})$ is a morphism of left $\mathbb{H}$-modules, then $(f,f)$ is a morphism in ${}_{(id_{H}, \Gamma_{H_{1}})}{\sf Mod}$ between $(M,M,\psi_{M}^{1},\Gamma_{M},\psi_{M}^{2},id_{M})$ and $(N,N,\psi_{N}^{1},\Gamma_{N},\psi_{N}^{2},id_{N})$, what follows by the fact that $f\colon (M,\Gamma_{M})\rightarrow (N,\Gamma_{N})$ is a morphism of left $H_{2}$-modules, which is straightforward to show.
\end{proof}

\begin{corollary}
Let $\left(T\begin{array}{c}H\\\downarrow\\B\end{array},\varphi_{H}\right)$ be a relative Rota-Baxter operator in ${\sf rRB}^{\star}$. If $T\colon H\rightarrow B$ is an isomorphism, then there exists a functor 
\[{\sf V}\colon {}_{\overline{\mathbb{H}}}{\sf Mod}\longrightarrow {}_{(T,\varphi_{H})}{\sf Mod}\]
where $\overline{\mathbb{H}}={\sf G}\left(\left(T\begin{array}{c}H\\\downarrow\\B\end{array},\varphi_{H}\right)\right)$ is the Hopf brace introduced in Theorem \ref{adjunction}.
\end{corollary}
\begin{proof}
Note that when $T$ is an isomorphism, then $(id_{H},T^{-1})$ is a morphism of relative Rota-Baxter operators between $\left(T\begin{array}{c}H\\\downarrow\\B\end{array},\varphi_{H}\right)$ and $\left(id_{H}\begin{array}{c}H\\\downarrow\\\overline{H}\end{array},\overline{\Gamma}_{H}\overset{\eqref{overlGammaH}}{=}\varphi_{H}\circ(T\otimes H)\right)$. Therefore, we can define {\sf V} as the following composition of functors:
\[{\sf V}\coloneqq {\sf R}_{(id_{H},T^{-1})}\circ{\sf W}_{\overline{\mathbb{H}}}\]
which is defined on objects by 
\[{\sf V}((M,\overline{\psi}_{M}^{1},\overline{\psi}_{M}^{2}))=(M,M,\overline{\psi}_{M}^{1},\overline{\Gamma}_{M}\circ(T^{-1}\otimes M),\overline{\psi}_{M}^{2}\circ(T^{-1}\otimes M),id_{M})\]
and on morphisms by ${\sf V}(f)=(f,f)$.
\end{proof}

\begin{theorem}
Let $\left(T\begin{array}{c}H\\\downarrow\\B\end{array},\varphi_{H}\right)$ be a relative Rota-Baxter operator in ${\sf rRB}^{\star}$. There exists a functor
\[{\sf U}\colon {}_{(T,\varphi_{H})}{\sf Mod}\longrightarrow {}_{\overline{\mathbb{H}}}{\sf Mod},\]
where $\overline{\mathbb{H}}={\sf G}\left(\left(T\begin{array}{c}H\\\downarrow\\B\end{array},\varphi_{H}\right)\right)$ is the Hopf brace introduced in Theorem \ref{adjunction}, defined on objects by $${\sf U}((M,N,\phi_{M},\varphi_{M},\varphi_{N}, \gamma))=(M,\phi_{M},\overline{\varphi}_{M})$$ being $\overline{\varphi}_{M}\coloneqq\phi_{M}\circ(H\otimes (\varphi_{M}\circ(T\otimes M)))\circ(\delta_{H}\otimes M)$, and on morphisms by ${\sf U}(r,s)=r$.
\end{theorem}
\begin{proof}
Consider $(M,N,\phi_{M},\varphi_{M},\varphi_{N}, \gamma)$ an object in ${}_{(T,\varphi_{H})}{\sf Mod}$ and let's show that $(M,\phi_{M},\overline{\varphi}_{M})$ is a module over the Hopf brace $\overline{\mathbb{H}}$. At first, we will see that $(M,\overline{\varphi}_{M})$ is a left $\overline{H}$-module. Indeed, it results straightforward to show that $\overline{\varphi}_{M}\circ(\eta_{H}\otimes M)=id_{M}$ and also
\begin{align*}
&\overline{\varphi}_{M}\circ(H\otimes\overline{\varphi}_{M})\\=&\phi_{M}\circ(H\otimes(\phi_{M}\circ ((\varphi_{H}\circ (T\otimes H))\otimes(\varphi_{M}\circ(T\otimes(\varphi_{M}\circ(T\otimes M)))))\circ (H\otimes c_{H,H}\otimes H\otimes M)\\&\circ(\delta_{H}\otimes H\otimes H\otimes M)))\circ(\delta_{H}\otimes\delta_{H}\otimes M)\;\footnotesize\textnormal{(by \eqref{compatmodHmodB}, the condition of coalgebra morphism for $T$ and naturality of $c$)}\\=&\phi_{M}\circ(\overline{\mu}_{H}\otimes (\varphi_{M}\circ((\mu_{B}\circ (T\otimes T))\otimes M)))\circ (((H\otimes c_{H,H}\otimes H)\circ(\delta_{H}\otimes\delta_{H}))\otimes M)\;\footnotesize\textnormal{(by the module axioms for}\\&\footnotesize\textnormal{$(M,\phi_{M})$ and $(M,\varphi_{M})$ and coassociativity of $\delta_{H}$)}\\=&\phi_{M}\circ(H\otimes(\varphi_{M}\circ(T\otimes M)))\circ(((\overline{\mu}_{H}\otimes\overline{\mu}_{H})\circ(H\otimes c_{H,H}\otimes H)\circ(\delta_{H}\otimes\delta_{H}))\otimes M)\;\footnotesize\textnormal{(by \eqref{rRBcond})}\\=&\overline{\varphi}_{M}\circ(\overline{\mu}_{H}\otimes M)\;\footnotesize\textnormal{(by the condition of coalgebra morphism for $\overline{\mu}_{H}$)}.
\end{align*}
	
So, to conclude that {\sf U} is well-defined on objects we have to see that the triple $(M,\phi_{M},\overline{\varphi}_{M})$ satisfies \eqref{compatmodbrace}. Note that, in this situation,
\begin{equation}\label{Gammabarra}
\overline{\Gamma}_{M}=\varphi_{M}\circ(T\otimes M)
\end{equation}
as we will prove in what follows:
\begin{align*}
&\overline{\Gamma}_{M}\\=&\phi_{M}\circ((\lambda_{H}\ast id_{H})\otimes (\varphi_{M}\circ(T\otimes M)))\circ(\delta_{H}\otimes M)\;\footnotesize\textnormal{(by coassociativity of $\delta_{H}$ and module axioms for $(M,\phi_{M})$)}\\=&\varphi_{M}\circ(T\otimes M)\;\footnotesize\textnormal{(by \eqref{antipode} and (co)unit properties)}.
\end{align*}
Then, \eqref{compatmodbrace} follows by
\begin{align*}
&\phi_{M}\circ(\overline{\mu}_{H}\otimes\overline{\Gamma}_{M})\circ(H\otimes c_{H,H}\otimes M)\circ(\delta_{H}\otimes H\otimes M)\\=&\phi_{M}\circ(H\otimes(\phi_{M}\circ((\varphi_{H}\circ(T\otimes H))\otimes(\varphi_{M}\circ(T\otimes M)))\circ (H\otimes c_{H,H}\otimes M)\circ(\delta_{H}\otimes H\otimes M)))\circ(\delta_{H}\otimes H\otimes M)\\&\footnotesize\textnormal{(by \eqref{Gammabarra}, module axioms for $(M,\phi_{M})$ and coassociativity of $\delta_{H}$)}\\=&\phi_{M}\circ(H\otimes(\phi_{M}\circ(\varphi_{H}\otimes\varphi_{M})\circ(B\otimes c_{B,H}\otimes M)\circ((\delta_{B}\circ T)\otimes H\otimes M)))\circ(\delta_{H}\otimes H\otimes M)\;\footnotesize\textnormal{(by naturality of $c$}\\&\footnotesize\textnormal{and the condition of coalgebra morphism for $T$)}\\=&\overline{\varphi}_{M}\circ(H\otimes\phi_{M})\;\footnotesize\textnormal{(by \eqref{compatmodHmodB})}.
\end{align*}
	
On the other hand, to show that {\sf U} is well-defined on morphisms it is enough to see that if $$(r,s)\colon (M,N,\phi_{M},\varphi_{M},\varphi_{N},\gamma)\rightarrow (P,Q,\phi_{P},\varphi_{P},\varphi_{Q},\theta)\in{}_{(T,\varphi_{H})}{\sf Mod},$$ then $r\colon (M,\overline{\varphi}_{M})\rightarrow (P,\overline{\varphi}_{P})$ is a morphism of left $\overline{H}$-modules. Indeed,
\begin{align*}
&r\circ\overline{\varphi}_{M}\\=&\phi_{P}\circ(H\otimes(r\circ\varphi_{M}\circ (T\otimes M)))\circ(\delta_{H}\otimes M)\;\footnotesize\textnormal{(by the condition of morphism of left $H$-modules}\\&\footnotesize\textnormal{for $f\colon (M,\phi_{M})\rightarrow (P,\phi_{P})$)}\\=&\overline{\varphi}_{P}\circ(H\otimes f)\;\footnotesize\textnormal{(by the condition of morphism of left $B$-modules for $r\colon (M,\varphi_{M})\rightarrow (P,\varphi_{P})$)}.\qedhere
\end{align*}
\end{proof}

\begin{theorem}\label{adjVU}
Let $\left(T\begin{array}{c}H\\\downarrow\\B\end{array},\varphi_{H}\right)$ be a relative Rota-Baxter operator in ${\sf rRB}^{\star}$ and assume that $T\colon H\rightarrow B$ is an isomorphism. Under these hypothesis, the functor {\sf V} is left adjoint of {\sf U}.
\end{theorem}
\begin{proof}
Let $(M,\overline{\psi}_{M}^{1},\overline{\psi}_{M}^{2})$ be a module over the Hopf brace $\overline{\mathbb{H}}$ and $(P,Q,\phi_{P},\varphi_{P},\varphi_{Q},\theta)$ a module over the relative Rota-Baxter operator $\left(T\begin{array}{c}H\\\downarrow\\B\end{array},\varphi_{H}\right)$. We have to set a bijection ${}^{(M,\overline{\psi}_{M}^{1},\overline{\psi}_{M}^{2})}\Lambda_{(P,Q,\theta)}$ between 
\[\Hom_{{}_{\overline{\mathbb{H}}}{\sf Mod}}((M,\overline{\psi}_{M}^{1},\overline{\psi}_{M}^{2}),(P,\phi_{P},\overline{\varphi}_{P}))\]
and
\[\Hom_{{}_{(T,\varphi_{H})}{\sf Mod}}((M,M,\overline{\psi}_{M}^{1},\overline{\Gamma}_{M}\circ(T^{-1}\otimes M),\overline{\psi}_{M}^{2}\circ(T^{-1}\otimes M),id_{M}),(P,Q,\phi_{P},\varphi_{P},\varphi_{Q},\theta)).\]
	
On the one hand, take $f\colon (M,\overline{\psi}_{M}^{1},\overline{\psi}_{M}^{2})\rightarrow (P,\phi_{P},\overline{\varphi}_{P})\in{}_{\overline{\mathbb{H}}}{\sf Mod}$ and let's see that $(f,\theta\circ f)$ is a morphism in ${}_{(T,\varphi_{H})}{\sf Mod}$ between $(M,M,\overline{\psi}_{M}^{1},\overline{\Gamma}_{M}\circ(T^{-1}\otimes M),\overline{\psi}_{M}^{2}\circ(T^{-1}\otimes M),id_{M})$ and $(P,Q,\phi_{P},\varphi_{P},\varphi_{Q},\theta)$. Then, we have to show that $f\colon (M,\overline{\Gamma}_{M}\circ(T^{-1}\otimes M))\rightarrow (P,\varphi_{P})$ and $\theta\circ f\colon (M,\overline{\psi}_{M}^{2}\circ(T^{-1}\otimes M))\rightarrow (Q,\varphi_{Q})$ are morphisms of left $B$-modules, what follows by
\begin{align*}
&f\circ\overline{\Gamma}_{M}\circ(T^{-1}\otimes M)\\=&\phi_{P}\circ (\lambda_{H}\otimes (f\circ\overline{\psi}_{M}^{2}))\circ((\delta_{H}\circ T^{-1})\otimes M)\;\footnotesize\textnormal{(by the condition of morphism of left $H$-modules for $f\colon (M,\overline{\psi}_{M}^{1})\rightarrow (P,\phi_{P})$)}\\=&\overline{\Gamma}_{P}\circ (T^{-1}\otimes f)\;\footnotesize\textnormal{(by the condition of morphism of left $\overline{H}$-modules for $f\colon (M,\overline{\psi}_{M}^{2})\rightarrow (P,\overline{\varphi}_{P})$)}\\=&\varphi_{P}\circ((T\circ T^{-1})\otimes f)\;\footnotesize\textnormal{(by \eqref{Gammabarra})}\\=&\varphi_{P}\circ(B\otimes f),
\end{align*}
and
\begin{align*}
&\theta\circ f\circ \overline{\psi}_{M}^{2}\circ (T^{-1}\otimes M)\\=&\theta\circ \overline{\varphi}_{P}\circ (T^{-1}\otimes f)\;\footnotesize\textnormal{(by the condition of morphism of left $\overline{H}$-modules for $f\colon (M,\overline{\psi}_{M}^{2})\rightarrow (P,\overline{\varphi}_{P})$)}\\=&\varphi_{Q}\circ ((T\circ T^{-1})\otimes(\theta\circ f))\;\footnotesize\textnormal{(by \eqref{compatmodT})}\\=&\varphi_{Q}\circ(B\otimes (\theta\circ f)).
\end{align*}
	
Therefore, let's define
\[{}^{(M,\overline{\psi}_{M}^{1},\overline{\psi}_{M}^{2})}\Lambda_{(P,Q,\theta)}(f)=(f,\theta\circ f).\]
	
On the other hand, consider $(r,s)\colon (M,M,\overline{\psi}_{M}^{1},\overline{\Gamma}_{M}\circ(T^{-1}\otimes M),\overline{\psi}_{M}^{2}\circ(T^{-1}\otimes M),id_{M})\rightarrow (P,Q,\phi_{P},\varphi_{P},\varphi_{Q},\theta)$ a morphism in ${}_{(T,\varphi_{H})}{\sf Mod}$ and let's prove that $r$ is a morphism in ${}_{\overline{\mathbb{H}}}{\sf Mod}$ between $(M,\overline{\psi}_{M}^{1},\overline{\psi}_{M}^{2})$ and $(P,\phi_{P},\overline{\varphi}_{P})$. To see this fact it is enough to compute that $r\colon (M,\overline{\psi}_{M}^{2})\rightarrow (P,\overline{\varphi}_{P})$ is a morphism of left $\overline{H}$-modules. Indeed,
\begin{align*}
&\overline{\varphi}_{P}\circ(H\otimes r)\\=&\phi_{P}\circ (H\otimes(r\circ\overline{\Gamma}_{M}))\circ(\delta_{H}\otimes M)\;\footnotesize\textnormal{(by the condition of morphism of left $B$-modules for $r\colon (M,\overline{\Gamma}_{M}\circ(T^{-1}\otimes M))\rightarrow (P,\varphi_{P})$)}\\=&r\circ \overline{\psi}_{M}^{1}\circ(H\otimes \overline{\Gamma}_{M})\circ(\delta_{H}\otimes M)\;\footnotesize\textnormal{(by the condition of morphism of left $H$-modules for $r\colon (M,\overline{\psi}_{M}^{1})\rightarrow (P,\phi_{P})$)}\\=&r\circ\overline{\psi}_{M}^{2}\;\footnotesize\textnormal{(by \eqref{psi2expression})}.
\end{align*}

Then, we define
\[({}^{(M,\overline{\psi}_{M}^{1},\overline{\psi}_{M}^{2})}\Lambda_{(P,Q,\theta)})^{-1}(r,s)=r.\]
	
So, ${}^{(M,\overline{\psi}_{M}^{1},\overline{\psi}_{M}^{2})}\Lambda_{(P,Q,\theta)}$ is a bijection because 
\begin{align*}
&(({}^{(M,\overline{\psi}_{M}^{1},\overline{\psi}_{M}^{2})}\Lambda_{(P,Q,\theta)})^{-1}\circ{}^{(M,\overline{\psi}_{M}^{1},\overline{\psi}_{M}^{2})}\Lambda_{(P,Q,\theta)})(f)=({}^{(M,\overline{\psi}_{M}^{1},\overline{\psi}_{M}^{2})}\Lambda_{(P,Q,\theta)})^{-1}(f,\theta\circ f)=f,
\end{align*}
and 
\begin{align*}
&({}^{(M,\overline{\psi}_{M}^{1},\overline{\psi}_{M}^{2})}\Lambda_{(P,Q,\theta)}\circ({}^{(M,\overline{\psi}_{M}^{1},\overline{\psi}_{M}^{2})}\Lambda_{(P,Q,\theta)})^{-1})(r,s)={}^{(M,\overline{\psi}_{M}^{1},\overline{\psi}_{M}^{2})}\Lambda_{(P,Q,\theta)}(r)=(r,\theta\circ r)\\=&(r,s)\;\footnotesize\textnormal{(by condition (iv) of morphism in ${}_{(T,\varphi_{H})}{\sf Mod}$)}.\qedhere
\end{align*}
\end{proof}

\begin{definition}
Let $\left(T\begin{array}{c}H\\\downarrow\\B\end{array},\varphi_{H}\right)$ be a relative Rota-Baxter operator. We will define ${}_{(T,\varphi_{H})}{\sf Mod^{iso}}$ as the full subcategory of ${}_{(T,\varphi_{H})}{\sf Mod}$ whose objects, $(M,N,\phi_{M},\varphi_{M},\varphi_{N},\gamma)$, satisfy that $\gamma\colon M\rightarrow N$ is an isomorphism in {\sf C}.	
\end{definition}

\begin{theorem}\label{equivVU}
Let $\left(T\begin{array}{c}H\\\downarrow\\B\end{array},\varphi_{H}\right)$ be a relative Rota-Baxter operator in ${\sf rRB}^{\star}$  such that $T\colon H\rightarrow B$ is an isomorphism. The categories ${}_{(T,\varphi_{H})}{\sf Mod^{iso}}$ and ${}_{\overline{\mathbb{H}}}{\sf Mod}$ are equivalent.
\end{theorem}
\begin{proof}
Define the functor ${\sf U}'$ as the restriction of ${\sf U}$ to ${}_{(T,\varphi_{H})}{\sf Mod^{iso}}$. It result straightforward to show that 
\[{\sf U'}\circ{\sf V}={\sf id}_{{}_{\overline{\mathbb{H}}}{\sf Mod}}.\]
	
On the other hand, consider $(M,N,\phi_{M},\varphi_{M},\varphi_{N},\gamma)$ in ${}_{(T,\varphi_{H})}{\sf Mod^{iso}}$. We have that
\begin{align*}
&({\sf V}\circ{\sf U}')((M,N,\phi_{M},\varphi_{M},\varphi_{N},\gamma))\\=&{\sf V}((M,\phi_{M},\overline{\varphi}_{M}))\\=&(M,M,\phi_{M},\varphi_{M},\overline{\varphi}_{M}\circ(T^{-1}\otimes M),id_{M})\;\footnotesize\textnormal{(by \eqref{Gammabarra})},
\end{align*}
which is isomorphic to $(M,N,\phi_{M},\varphi_{M},\varphi_{N},\gamma)$ via $(id_{M},\gamma^{-1})$. Indeed, to show that $(id_{M},\gamma^{-1})$ is a morphism in ${}_{(T,\varphi_{H})}{\sf Mod}$ it is enough to see that $\gamma^{-1}\colon (N,\varphi_{N})\rightarrow (M,\overline{\varphi}_{M}\circ (T^{-1}\otimes M))$ is a morphism of left $B$-modules, what follows by
\begin{align*}
&\gamma^{-1}\circ\varphi_{N}\\=&\gamma^{-1}\circ\varphi_{N}\circ(T\otimes\gamma)\circ(T^{-1}\otimes\gamma^{-1})\;\footnotesize\textnormal{(by the fact that $T$ and $\gamma$ are isomorphisms)}\\=&\overline{\varphi}_{M}\circ (T^{-1}\otimes\gamma^{-1})\;\footnotesize\textnormal{(by \eqref{compatmodT})}.
\end{align*}
As a consequence,
\[{\sf V}\circ{\sf U}'\simeq {\sf id}_{{}_{(T,\varphi_{H})}{\sf Mod^{iso}}}\]
what concludes the proof.
\end{proof}

\begin{corollary} 
Let $\mathbb{H}$ be a cocommutative Hopf brace. The category ${}_{(id_{H},\Gamma_{H_{1}})}{\sf Mod^{iso}}$ is equivalent to ${}_{\mathbb{H}}{\sf Mod}$.
\end{corollary}
\begin{proof}
The proof is a direct consequence of previous theorem taking into account that $({\sf G}\circ {\sf F})(\mathbb{H})=\mathbb{H}$, where ${\sf F}$ and $\sf{G}$ are the functors introduced in Theorem \ref{adjunction}.
\end{proof}

\section{Projections of relative Rota-Baxter operators}\label{sect3}

This section is devoted to the study of projections between relative Rota-Baxter operators. A projection between this kind of objects involves projections of Hopf algebras. Then, at the beginning of this section we will make a brief summary of the basic theory linked to projections of Hopf algebras including the Yetter-Drinfeld modules following  \cite{MN}, \cite{MAJ2} and \cite{RAD}.

\begin{definition}
{\rm 
Let $X$ be a Hopf algebra in ${\sf  C}$. We shall denote by $^{\sf X}_{\sf X}{\sf  Y}{\sf  D}$ the category of left Yetter-Dinfeld modules over $X$. More concretely, a triple $M=(M, \varphi_{M}, \rho_{M})$ is an object in $^{\sf X}_{\sf X}{\sf  Y}{\sf  D}$ if $(M,\varphi_{M})$ is a left $X$-module, $(M,\rho_{M})$ is a left $X$-comodule and the following identity
\begin{gather}
(\mu_{X}\otimes M)\circ (X\otimes c_{M,X})\circ ((\rho_{M}\circ \varphi_{M})\otimes X)\circ (X\otimes c_{X,M})\circ (\delta_{X}\otimes M)\\\nonumber=(\mu_{X}\otimes \varphi_{M})\circ (X\otimes c_{X,X}\otimes M)\circ (\delta_{X}\otimes \rho_{M}).\end{gather}
holds. The morphisms in $\;^{\sf X}_{\sf X}{\sf  Y}{\sf  D}$ are morphisms of left $X$-modules and left $X$-comodules.

For example, for any Hopf algebra $X$, $(X, \varphi_{X}^{ad}, \rho_{X}=\delta_{X})$ and $(X, \varphi_{X}=\mu_{X}, \rho_{X}^{ad})$ are left Yetter-Drinfeld modules over $X$. Also, any left $X$-module $(M,\varphi_{M})$ over a cocommutative Hopf algebra $X$ is a Yetter-Drinfeld module with the trivial left coaction $\rho_{M}=\eta_{X}\ot M$. Finally, the triple $(M, \varphi_{M}=\varepsilon_{X}\ot M, \rho_{M}=\eta_{X}\ot M)$ is a left Yetter-Drinfeld module for all Hopf algebra $X$. 
	}
\end{definition}

The category $^{\sf X}_{\sf X}{\sf  Y}{\sf  D}$ is strict monoidal with the usual tensor product in ${\sf  C}$, that is to say, for
$M$, $N$ in ${\sf {}^{X}_{X}YD}$, $M\otimes N$ is a left Yetter-Drinfeld module over $X$ with the tensor module and comodule structures given by \begin{gather*}\varphi_{M\otimes N}=(\varphi_{M}\otimes \varphi_{N})\circ (X\otimes c_{X,M}\otimes N)\circ (\delta_{X}\otimes M\otimes N),
\\\rho_{M\otimes N}=(\mu_{X}\otimes M\otimes N)\circ (X\otimes c_{M,X}\otimes N)\circ (\rho_{M}\otimes \rho_{N}).\end{gather*}

If the antipode of $X$ is an isomorphism,  $\;^{\sf X}_{\sf X}{\sf  Y}{\sf  D}$ is a braided monoidal category category where the braiding $t_{M,N}:M\otimes N\rightarrow N\otimes N$ is given by $t_{M,N}=(\varphi_{N}\otimes M)\circ (X\otimes c_{M,N})\circ (\rho_{M}\otimes N).$ It is immediate to see that $t_{M,N}$ is natural and it is an isomorphism with inverse
$$t_{M,N}^{-1}=c_{M,N}^{-1}\circ (\varphi_{N}\otimes M)\circ (\lambda_{X}^{-1}\otimes N\otimes M)\circ (c_{X,N}^{-1}\otimes M)\circ (N\otimes \rho_{M}).$$

\begin{definition}
A projection of Hopf algebras in {\sf C} is a 4-tuple $(X,Y, f,g)$ where $X$, $Y$ are Hopf algebras, and $f:X\rightarrow Y$, $g:Y\rightarrow X$ are Hopf algebra morphisms such that $g\co f=id_{X}$. 
	
A morphism between  projections of Hopf algebras $(X,Y, f,g)$ and $(X^{\prime},Y^{\prime}, f^{\prime},g^{\prime})$ is a pair $(x,y)$, where $x:X\rightarrow X^{\prime}$, $y:Y\rightarrow Y^{\prime}$ are Hopf algebra morphisms such that
\begin{equation}
\label{proj-morph}
y\circ f=f^{\prime}\circ x, \;\; x\circ g=g^{\prime}\circ y. 
\end{equation}
	
With the obvious composition of morphisms we can define a category whose objects are  Hopf algebra projections and whose morphisms are morphisms of Hopf algebra projections. We denote this category by ${\sf  P}({\sf  Hopf})$.
	
It is obvious that there exists a functor ${\sf P}_{triv}:{\sf Hopf}\rightarrow {\sf P(Hopf)}$ defined on objects by ${\sf P}_{triv}(X)=(X,X, id_{X}, id_{X})$ and on morphisms by ${\sf P}(f)=(f,f)$.
\end{definition}

Let $(X, Y, f, g)$ be an object in ${\sf  P}({\sf  Hopf})$.   The morphism 
$q_{Y}\coloneqq id_{Y}\ast (f\co \lambda_{X}\co g)$ is  idempotent and, as a consequence, there exists an object $I(q_{Y})$, called the object of coinvariants, an epimorphism $p_{Y}$ and a monomorphism $i_{Y}$ such that $q_{Y}=i_{Y}\circ p_{Y}$ and $p_{Y}\circ i_{Y}=id_{I(q_{Y})}$. As a consequence,
$$
\xymatrix{&\IY\ar[rr]^-{i_{Y}} & &Y\ar@<1ex>[rr]^-{(Y\otimes g)\circ \delta_{Y}}\ar@<-1ex>[rr]_-{Y\otimes\eta_{X}} & &Y\otimes X}
$$
is an equalizer diagram and $I(q_{Y})$ is a left $X$-module algebra where  the algebra structure is defined by
\begin{equation}
\label{alg-id} 
\eta_{I(q_{Y})}=p_{Y}\co \eta_{Y},\;\;\; \mu_{I(q_{Y})}=p_{Y}\co \mu_{Y}\co (i_{Y}\ot i_{Y}),
\end{equation}
i.e., $\eta_{I(q_{Y})}$ is the unique morphism such that $i_{Y}\co \eta_{I(q_{Y})}=\eta_{Y}$ and $\mu_{I(q_{Y})}$ is the unique morphism such that 
\begin{equation}
\label{m-i}
i_{Y}\co \mu_{I(q_{Y})}=\mu_{Y}\co (i_{Y}\ot i_{Y}).
\end{equation} 
The action $\psi_{I(q_{Y})}: X\ot I(q_{Y})\rightarrow I(q_{Y})$ is 
$
\psi_{I(q_{Y})}=p_{Y}\co\mu_{Y}\co (f\ot i_{Y}),
$
and then $\psi_{I(q_{Y})}$ is the unique morphism such that 
\begin{equation}
\label{phi2}
i_{Y}\co \psi_{I(q_{Y})}=\varphi_{Y}^{ad}\co (f\ot i_{Y}).
\end{equation}

On the other hand,
$$
\xymatrix{&Y\otimes X\ar@<1ex>[rr]^-{\mu_{Y}\circ(Y\otimes f)}\ar@<-1ex>[rr]_-{Y\otimes\varepsilon_{X}} & &Y\ar[rr]^-{p_{Y}} & &\IY}
$$
is a coequalizer diagram and, as a consequence, $I(q_{Y})$ is a left $X$-comodule coalgebra with
\begin{equation}
\label{coalg-id} 
\varepsilon_{I(q_{Y})}=\varepsilon_{Y}\co i_{Y},\;\; \delta_{I(q_{Y})}=(p_{Y}\otimes p_{Y})\circ \delta_{Y}\co i_{Y}
\end{equation}
and coaction $\rho_{I(q_{Y})}: I(q_{Y})\rightarrow X\ot I(q_{Y})$ defined by 
$$
\rho_{I(q_{Y})}=(g\otimes p_{Y})\circ \delta_{Y}\circ i_{Y}.
$$
In this case $\varepsilon_{I(q_{Y})}$ is the unique morphism such that $\varepsilon_{I(q_{Y})}\co p_{Y}=\varepsilon_{Y}$, $\delta_{I(q_{Y})}$ is the unique morphism such that 
\begin{equation}
\label{co-i}
\delta_{I(q_{Y})}\co p_{Y}=(p_{Y}\ot p_{Y})\co \delta_{Y},
\end{equation}
and  the coaction $\rho_{I(q_{Y})}$ is the unique morphism satisfying 
\begin{equation}
\rho_{I(q_{Y})}\co p_{Y}=(g\ot p_{Y})\co \rho^{ad}_{Y}.
\end{equation}

The algebra-coalgebra $I(q_{Y})$, with the action $\psi_{I(q_{Y})}$ and the coaction $\rho_{I(q_{Y})}$, is a Hopf algebra in   $^{\sf X}_{\sf X}{\sf  Y}{\sf  D}$ with antipode $$\lambda_{I(q_{Y})}=\psi_{I(q_{Y})}\circ (X\otimes (p_{Y}\circ \lambda_{Y}\circ i_{Y}))\circ \rho_{I(q_{Y})}.$$

Also, using that $i_{Y}$ is an equalizer morphism and $p_{Y}$ is a coequalizer, we obtain the following identities:
\begin{equation}
\label{id1}
p_{Y}\co \mu_{Y}\co (Y\ot q_{Y})=p_{Y}\co \mu_{Y}, \;\;\;\;\;
(Y\ot q_{Y})\co \delta_{Y}\co i_{Y}= \delta_{Y}\co i_{Y}. 
\end{equation}

Note that $i_{Y}$ is a coalgebra morphism iff 
\begin{equation}
\label{iz-coal}
(q_{Y}\ot Y)\co \delta_{Y}\co i_{Y}=\delta_{Y}\co i_{Y}.
\end{equation}

If $Y$ is cocommutative, condition (\ref{iz-coal}) always holds. This fact was proved by Sweedler in \cite{SW} for projections of Hopf algebras in a category of vector spaces. On the other hand, there exist examples where  $i_{Y}$ is not a coalgebra morphism (see \cite{BCM} for the complete details). In any case, if $i_Y$ is a coalgebra morphism, then we have that $I(q_{Y})$ is a Hopf algebra in ${\sf C}$ because $\rho_{I(q_{Y})}$ is trivial.

Similarly, $p_{Y}$ is an algebra morphism iff 
\begin{equation}
\label{pz-al}
p_{Y}\co \mu_{Y}\co (q_{Y}\ot Y)=p_{Y}\co \mu_{Y}.
\end{equation}

Equivalently, $p_{Y}$ is an algebra morphism iff  $\psi_{I(q_{Y})}=\varepsilon_{X}\ot I(q_{Y})$ (see \cite{MN}). Therefore, if $p_Y$ is an algebra morphism, then we have that $I(q_{Y})$ is a Hopf algebra in ${\sf C}$ because $\psi_{I(q_{Y})}$ is trivial. 

Finally, if $(X, Y, f, g)$ is in ${\sf  P}({\sf  Hopf})$ and $Y$ is cocommutative, then the morphism $q_{Y}$ is a coalgebra morphism. Also, under these conditions, $i_{Y}$ is a coalgebra morphism and the following equality 
\begin{equation}
\label{lxy}
i_{Y}\circ \lambda_{I(q_{Y})}=\lambda_{Y}\circ i_{Y}.
\end{equation}
holds (see \cite{FGRR}).

In the following definition the notion of projection between Hopf algebras is extended in order to introduce what a projection between relative Rota-Baxter operators is.
\begin{definition}\label{defprojrB}
Let $\left(T\begin{array}{c}H\\\downarrow\\B \end{array},\varphi_{H}\right)$ and $\left(L\begin{array}{c}A\\\downarrow\\D \end{array},\varphi_{A}\right)$ be relative Rota-Baxter operators. We will say that a 6-tuple 
\[\left(\left(T\begin{array}{c}H\\\downarrow\\B \end{array},\varphi_{H}\right),\left(L\begin{array}{c}A\\\downarrow\\D \end{array},\varphi_{A}\right),f, h,g,l\right)\]
is a projection of relative Rota-Baxter operators if
\begin{itemize}
\item[(i)] The 4-tupla $(H,A,f,g)$ is an object in ${\sf P(Hopf)}$,
\item[(ii)] The 4-tupla $(B,D,h,l)$ is an object in ${\sf P(Hopf)}$,
\item[(iii)] The pair $$(f,h)\colon \left(T\begin{array}{c}H\\\downarrow\\B \end{array},\varphi_{H}\right)\rightarrow \left(L\begin{array}{c}A\\\downarrow\\D \end{array},\varphi_{A}\right)$$ is a morphism of relative Rota-Baxter operators,
\item[(iv)] The pair $$(g,l)\colon \left(L\begin{array}{c}A\\\downarrow\\D \end{array},\varphi_{A}\right)\rightarrow \left(T\begin{array}{c}H\\\downarrow\\B \end{array},\varphi_{H}\right)$$ is a morphism of relative Rota-Baxter operators.
\end{itemize}

\begin{figure}[h]
\begin{tikzcd}
H \arrow[dd, "T"'] \arrow[rrrr, "f", shift left] \arrow["id_{H}", loop, distance=2em, in=145, out=215] &    &  &  & A \arrow[llll, "g", shift left] \arrow[dd, "L"] \\
{} \arrow[r, phantom, shift right=2]                                                                   & {} &  &  &                                                 \\
B \arrow[rrrr, "h", shift left] \arrow["id_{B}", loop, distance=2em, in=145, out=215]                  &    &  &  & D \arrow[llll, "l", shift left]                
\end{tikzcd}
\caption{Projection of relative Rota-Baxter operators.}
\end{figure}

Projections of relative Rota-Baxter operators give rise to a category whose morphisms are defined as follows: We will say that a 4-tuple $(x,y,z,t)$ is a morphism between the projections of relative Rota-Baxter operators
\[\left(\left(T\begin{array}{c}H\\\downarrow\\B \end{array},\varphi_{H}\right),\left(L\begin{array}{c}A\\\downarrow\\D \end{array},\varphi_{A}\right),f,h,g,l \right)\rightarrow \left(\left(T'\begin{array}{c}H'\\\downarrow\\B' \end{array},\varphi_{H'}\right),\left(L'\begin{array}{c}A'\\\downarrow\\D' \end{array},\varphi_{A'}\right),f',h',g',l' \right)\]
if the conditions
\begin{itemize}
\item[(i)] The pair $(x,y)\colon (H,A,f,g)\rightarrow (H',A',f',g')$ is a morphism in ${\sf P(Hopf)}$,
\item[(ii)] The pair $(z,t)\colon (B,D,h,l)\rightarrow (B',D',h',l')$ is a morphism in ${\sf P(Hopf)}$,
\item[(iii)] The pair $(x,z)\colon \left(T\begin{array}{c}H\\\downarrow\\B \end{array},\varphi_{H}\right)\rightarrow \left(T'\begin{array}{c}H'\\\downarrow\\B' \end{array},\varphi_{H'}\right)$ is a morphism in ${\sf rRB}$,
\item[(iv)] The pair $(y,t)\colon \left(L\begin{array}{c}A\\\downarrow\\D \end{array},\varphi_{A}\right)\rightarrow \left(L'\begin{array}{c}A'\\\downarrow\\D' \end{array},\varphi_{A'}\right)$ is a morphism in ${\sf rRB}$,
\end{itemize}
hold. This category will be denoted by ${\sf P(rRB)}$.

\begin{figure}[h]
\begin{tikzcd}
                                                                                               &  &  & H \arrow[dd, "T"'] \arrow[lllddd, "x"'] \arrow[rrrr, "f", shift left] \arrow["id_{H}", loop, distance=2em, in=145, out=215] &                                                      &  &  & A \arrow[llll, "g", shift left] \arrow[dd, "L"] \arrow[lllddd, "y"] \\
                                                                                               &  &  & {} \arrow[r, phantom, shift right=2]                                                                                        & {}                                                   &  &  &                                                                     \\
                                                                                               &  &  & B \arrow[lllddd, "z"'] \arrow[rrrr, "h", shift left] \arrow["id_{B}", loop, distance=2em, in=145, out=215]                  &                                                      &  &  & D \arrow[llll, "l", shift left] \arrow[lllddd, "t"]                 \\
H' \arrow[dd, "T'"'] \arrow[rrrr, "f'"] \arrow["id_{H'}", loop, distance=2em, in=145, out=215] &  &  &                                                                                                                             & A' \arrow[dd, "L'"] \arrow[llll, "g'", shift left=3] &  &  &                                                                     \\
                                                                                               &  &  &                                                                                                                             &                                                      &  &  &                                                                     \\
B' \arrow[rrrr, "h'", shift left] \arrow["id_{B}", loop, distance=2em, in=145, out=215]        &  &  &                                                                                                                             & D' \arrow[llll, "l'", shift left]                    &  &  &                                                                    
\end{tikzcd}
\caption{Morphism of projections of relative Rota-Baxter operators.}
\end{figure}
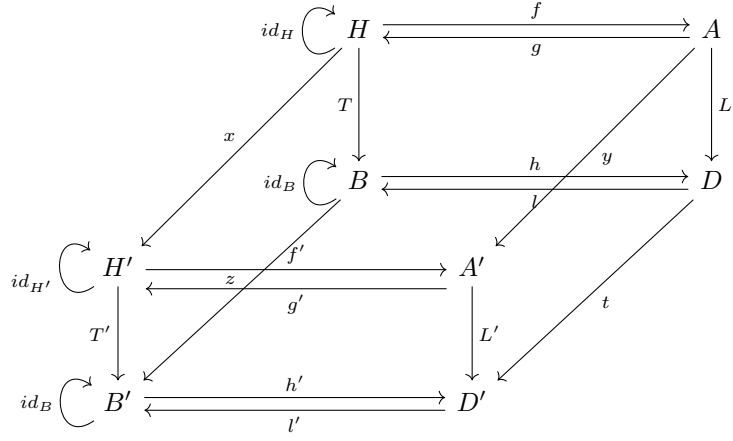

By ${\sf P(rRB^{\star})}$ we will denote to the full subcategory of ${\sf P(rRB)}$ where the involved relative Rota-Baxter operators  are objects in ${\sf rRB^{\star}}$. Finally ${\sf P({\sf coc}\textnormal{-}{\sf rRB})}$ denotes the full subcategory of ${\sf P(rRB^{\star})}$  where the involved objects are cocommutative relative Rota-Baxter operators.

\end{definition}

Let $$\left(\left(T\begin{array}{c}H\\\downarrow\\B \end{array},\varphi_{H}\right),\left(L\begin{array}{c}A\\\downarrow\\D \end{array},\varphi_{A}\right),f,h,g,l \right)$$ be an object in ${\sf P({\sf coc}\textnormal{-}{\sf rRB})}$. Due to being $(H,A,f,g)$ and $(B,D,h,l)$ objects in ${\sf P(Hopf)}$ with $A$ and $D$ cocommutative, we know that the respective objects of coinvariants, $\IA$ and $\ID$, are Hopf algebras in ${\sf C}$ and also the equalizers, $i_{A}$ and $i_{D}$, are coalgebra morphisms. Moreover, the morphism $L\circ i_{A}$ factors through the equalizer $i_{D}$ because 
\begin{align*}
&(D\otimes l)\circ\delta_{D}\circ L\circ i_{A}\\=&(L\otimes(l\circ L))\circ\delta_{A}\circ i_{A}\;\footnotesize\textnormal{(by the condition of coalgebra morphism for $L$)}\\=&(L\otimes T)\circ(A\otimes g)\circ\delta_{A}\circ i_{A}\;\footnotesize\textnormal{(by \eqref{cond1morrRB} for the morphism $(g,l)$)}\\=&((L\circ i_{A})\otimes(T\circ\eta_{H}))\;\footnotesize\textnormal{(by the equalizer condition for $i_{A}$)}\\=&(T\circ i_{A})\otimes\eta_{B}\;\footnotesize\textnormal{(by \eqref{etaT} for $T$)}.
\end{align*}
As a conclusion, there exists an unique $L_{0}\colon\IA\rightarrow\ID$ satisfying that \begin{equation}\label{T0cons} i_{D}\circ L_{0}=L\circ i_{A}. \end{equation} Composing on the left with the epimorphism $p_{D}$ and taking into account that $p_{D}\circ i_{D}=id_{\ID}$, we obtain that
\begin{equation}\label{T0def}
L_{0}=p_{D}\circ L\circ i_{A}.
\end{equation}
Also note that there exists an action $$\varphi_{\IA}\colon \ID\otimes\IA\rightarrow \IA,$$ which is obtained by factorization through the equalizer $i_{A}$ of the morphism $\varphi_{A}\circ(i_{D}\otimes i_{A})$. Indeed,
\begin{align*}
&(A\otimes g)\circ\delta_{A}\circ\varphi_{A}\circ(i_{D}\otimes i_{A})\\=&(A\otimes g)\circ(\varphi_{A}\otimes\varphi_{A})\circ(D\otimes c_{D,A}\otimes A)\circ((\delta_{D}\circ i_{D})\otimes(\delta_{A}\circ i_{A}))\;\footnotesize\textnormal{(by the condition of coalgebra morphism for $\varphi_{A}$)}\\=&(\varphi_{A}\otimes\varphi_{H})\circ(D\otimes c_{B,A}\otimes H)\circ (((D\otimes l)\circ\delta_{D}\circ i_{D})\otimes((A\otimes g)\circ\delta_{A}\circ i_{A}))\;\footnotesize\textnormal{(by \eqref{cond2morrRB} for $(g,l)$ and naturality of $c$)}\\=&(\varphi_{A}\circ(i_{D}\otimes i_{A}))\otimes\eta_{H}\;\footnotesize\textnormal{(by equalizer condition for $i_{D}$ and $i_{A}$, naturality of $c$ and module axioms)}.
\end{align*}
Therefore, there exists an unique $\varphi_{\IA}\colon \ID\otimes\IA\rightarrow \IA$ satisfying that
\begin{equation}\label{PhiIBcons}
i_{A}\circ\varphi_{\IA}=\varphi_{A}\circ(i_{D}\otimes i_{A}).
\end{equation}
Then, composing on the left with $p_{A}$ and using that $p_{A}\circ i_{A}=id_{\IA}$, the equality 
\begin{equation}\label{PhiIBdef}
\varphi_{\IA}=p_{A}\circ\varphi_{A}\circ(i_{D}\otimes i_{A})
\end{equation}
holds.

\begin{figure}[h]
\begin{minipage}{0.4  \linewidth}
\[\hspace{-3cm}\xymatrix{&\IA\ar@{-->}[d]^-{\exists^{\bullet}\,L_{0}}\ar[rr]^-{i_{A}} & &A\ar[d]^-{L}\ar@<1ex>[rr]^-{(A\otimes g)\circ\delta_{A}}\ar@<-1ex>[rr]_-{A\otimes\eta_{H}} & &A\otimes H\\&\ID\ar[rr]^-{i_{D}} & &D\ar@<1ex>[rr]^-{(D\otimes l)\circ\delta_{D}}\ar@<-1ex>[rr]_-{D\otimes\eta_{B}} & &D\otimes B}\]
\end{minipage}
\begin{minipage}{0.4 \linewidth}
\[\hspace{-1.5 cm}\xymatrix{&\IA\ar[rr]^-{i_{A}} & &A\ar@<1ex>[rr]^-{(A\otimes g)\circ\delta_{A}}\ar@<-1ex>[rr]_-{A\otimes\eta_{H}} & &A\otimes H\\& & &\ID\otimes\IA\ar@{-->}[llu]^-{\exists^{\bullet}\varphi_{\IA}}\ar[u]_-{\varphi_{A}\circ(i_{D}\otimes i_{A})} & &}\vspace{0.4cm}\]
\end{minipage}
\caption{The construction of $L_{0}$ and $\varphi_{\IA}$.}
\end{figure}

\begin{theorem}\label{thT0}
Let $$\left(\left(T\begin{array}{c}H\\\downarrow\\B \end{array},\varphi_{H}\right),\left(L\begin{array}{c}A\\\downarrow\\D \end{array},\varphi_{A}\right),f,h,g,l \right)$$ be an object in ${\sf P({\sf coc}\textnormal{-}{\sf rRB})}$. Then, if $L_{0}$ is the morphism defined by \eqref{T0def} and $\varphi_{\IA}$ the action introduced in \eqref{PhiIBdef}, we have that  $$\left(L_{0}\begin{array}{c}\IA\\\downarrow\\\ID \end{array},\varphi_{\IA}\right)$$ is a relative Rota-Baxter operator.
\end{theorem}

\begin{proof}
First of all note that $(\IA,\varphi_{\IA})$ is a left $\ID$-module algebra-coalgebra. Indeed, module axioms are straightforward thanks to \eqref{PhiIBcons} and module axioms for $(A,\varphi_{A})$. Moreover, it is direct to compute that $\eta_{\IA}$ is a morphism of left $\ID$-modules and
\begin{align*}
&i_{A}\circ\varphi_{\IA}\circ(\ID\otimes \mu_{\IA})\\=&\varphi_{A}\circ(D\otimes\mu_{A})\circ(i_{D}\otimes i_{A}\otimes i_{A})\;\footnotesize\textnormal{(by \eqref{PhiIBcons} and \eqref{m-i})}\\=&\mu_{A}\circ(\varphi_{A}\otimes\varphi_{A})\circ (D\otimes c_{D,A}\otimes A)\circ((\delta_{D}\circ i_{D})\otimes i_{A}\otimes i_{A})\;\footnotesize\textnormal{(by the condition of morphism of left $D$-modules for $\mu_{A}$)}\\=&\mu_{A}\circ((\varphi_{A}\circ(i_{D}\otimes i_{A}))\otimes(\varphi_{A}\circ(i_{D}\otimes i_{A})))\circ(\ID\otimes c_{\ID,\IA}\otimes \IA)\circ(\delta_{\ID}\otimes \IA\otimes\IA)\\&\footnotesize\textnormal{(by the condition of coalgebra morphism for $i_{D}$ and naturality of $c$)}\\=&i_{A}\circ \mu_{\IA}\circ(\varphi_{\IA}\otimes\varphi_{\IA})\circ(\ID\otimes c_{\ID,\IA}\otimes\IA)\circ(\delta_{\ID}\otimes\IA\otimes \IA)\;\footnotesize\textnormal{(by \eqref{PhiIBcons} and \eqref{m-i})},
\end{align*}
so $(\IA,\varphi_{\IA})$ is a left $\ID$-module algebra. To finish, $\varphi_{\IA}$ is a coalgebra morphism because 
\begin{align*}
&\delta_{\IA}\circ\varphi_{\IA}\\=&(p_{A}\otimes p_{A})\circ\delta_{A}\circ \varphi_{A}\circ(i_{D}\otimes i_{A})\;\footnotesize\textnormal{(by \eqref{PhiIBcons})}\\=&(p_{A}\otimes p_{A})\circ(\varphi_{A}\otimes\varphi_{A})\circ(D\otimes c_{D,A}\otimes A)\circ((\delta_{D}\circ i_{D})\otimes(\delta_{A}\circ i_{A}))\;\footnotesize\textnormal{(by the condition of coalgebra morphism for $\varphi_{A}$)}\\=&((p_{A}\circ\varphi_{A}\circ(i_{D}\otimes i_{A}))\otimes(p_{A}\circ\varphi_{A}\circ(i_{D}\otimes i_{A})))\circ(\ID\otimes c_{\ID,\IA}\otimes\IA)\circ(\delta_{\ID}\otimes\delta_{\IA})\\&\footnotesize\textnormal{(by the condition of coalgebra morphism for $i_{D}$ and $i_{A}$ and naturality of $c$)}\\=&(\varphi_{\IA}\otimes\varphi_{\IA})\circ(\ID\otimes c_{\ID,\IA}\otimes\IA)\circ(\delta_{\ID}\otimes\delta_{\IA})\;\footnotesize\textnormal{(by \eqref{PhiIBdef})},
\end{align*}
and 
\begin{align*}
&\varepsilon_{\IA}\circ\varphi_{\IA}\\=&\varepsilon_{A}\circ \varphi_{A}\circ (i_{D}\otimes i_{A})\;\footnotesize\textnormal{(by \eqref{coalg-id} and \eqref{PhiIBcons})}\\=&((\varepsilon_{D}\circ i_{D})\otimes(\varepsilon_{A}\circ i_{A}))\;\footnotesize\textnormal{(by the condition of coalgebra morphism for $\varphi_{A}$)}\\=&\varepsilon_{\ID}\otimes\varepsilon_{\IA}\;\footnotesize\textnormal{(by \eqref{coalg-id})}.
\end{align*}
	
On the other hand, let's see that $L_{0}$ is a coalgebra morphism. Indeed,
\begin{align*}
&\delta_{\ID}\circ L_{0}\\=&(p_{D}\otimes p_{D})\circ \delta_{D}\circ L\circ i_{A}\;\footnotesize\textnormal{(by \eqref{T0cons})}\\=&((p_{D}\circ L)\otimes(p_{D}\circ L))\circ\delta_{A}\circ i_{A}\;\footnotesize\textnormal{(by the condition of coalgebra morphism for $L$)}\\=&((p_{D}\circ L\circ i_{A})\otimes(p_{D}\circ L\circ i_{A}))\circ\delta_{\IA}\;\footnotesize\textnormal{(by the condition of coalgebra morphism for $i_{A}$)}\\=&(L_{0}\otimes L_{0})\circ\delta_{\IA}\;\footnotesize\textnormal{(by \eqref{T0def})}.
\end{align*}
	
So, the only condition that remains us to compute is \eqref{rRBcond}, which follows by:
\begin{align*}
&i_{D}\circ L_{0}\circ \mu_{\IA}\circ(\IA\otimes(\varphi_{\IA}\circ(L_{0}\otimes\IA)))\circ(\delta_{\IA}\otimes\IA)\\=&L\circ i_{A}\circ \mu_{\IA}\circ(\IA\otimes(\varphi_{\IA}\circ(L_{0}\otimes\IA)))\circ(\delta_{\IA}\otimes\IA)\;\footnotesize\textnormal{(by \eqref{T0cons})}\\=&L\circ \mu_{A}\circ (i_{A}\otimes(i_{A}\circ\varphi_{\IA}\circ(L_{0}\otimes\IA)))\circ(\delta_{\IA}\otimes\IA)\;\footnotesize\textnormal{(by \eqref{m-i})}\\=&L\circ\mu_{A}\circ (A\otimes (\varphi_{A}\circ(L\otimes A)))\circ(((i_{A}\otimes i_{A})\circ\delta_{\IA})\otimes i_{A})\;\footnotesize\textnormal{(by \eqref{PhiIBcons} and \eqref{T0cons})}\\=&L\circ\mu_{A}\circ (A\otimes (\varphi_{A}\circ(L\otimes A)))\circ(\delta_{A}\otimes A)\circ(i_{A}\otimes i_{A})\;\footnotesize\textnormal{(by the condition of coalgebra morphism for $i_{A}$)}\\=&\mu_{D}\circ((L\circ i_{A})\otimes(L\circ i_{A}))\;\footnotesize\textnormal{(by \eqref{rRBcond} for $\left(L\begin{array}{c}A\\\downarrow\\ D\end{array},\varphi_{A}\right)$)}\\=&i_{D}\circ\mu_{\ID}\circ(L_{0}\otimes L_{0})\;\footnotesize\textnormal{(by \eqref{T0cons} and \eqref{m-i})}.\qedhere
\end{align*}
\end{proof}

\begin{corollary}
Under the conditions of the previous theorem, the pair $(i_{A},i_{D})$ is a morphism of relative Rota-Baxter operators between $\left(L_{0}\begin{array}{c}\IA\\\downarrow\\\ID\end{array},\varphi_{\IA}\right)$ and $\left(L\begin{array}{c}A\\\downarrow\\D\end{array},\varphi_{A}\right)$.
\end{corollary}
\begin{proof}
By \eqref{T0cons}, \eqref{cond1morrRB} holds and \eqref{cond2morrRB} follows by \eqref{PhiIBcons}.
\end{proof}

\begin{example}
In \cite{Goncharov}, Goncharov proved that if $X$ is a cocommutative Hopf algebra, then $$\left(\lambda_{X}\begin{array}{c}X\\\downarrow\\X\end{array},\varphi_{X}^{ad}\right)$$ is a relative Rota-Baxter operator. So, if $H$ and $A$ are cocommutative Hopf algebras and $(H,A,f,g)\in{\sf P(Hopf)}$, then 
\[\left(\left(\lambda_{H}\begin{array}{c}H\\\downarrow\\H\end{array},\varphi_{H}^{ad}\right),\left(\lambda_{A}\begin{array}{c}A\\\downarrow\\A\end{array},\varphi_{A}^{ad}\right),f,f,g,g\right)\]
is a projection of relative Rota-Baxter operators. As a consequence, by Theorem \ref{thT0}, there exists a relative Rota-Baxter operator, $\left(L_{0}\begin{array}{c}\IA\\\downarrow\\\IA\end{array},\varphi_{\IA}\right)$, where $L_{0}$ and $\varphi_{\IA}$ satisfy \eqref{T0def} and \eqref{PhiIBdef}, respectively. Note that, due to being $A$ cocommutative, by (\ref{lxy}), 
\[L_{0}=p_{A}\circ\lambda_{A}\circ i_{A}=\lambda_{\IA}.\]
Then, we obtain that $\left(\lambda_{\IA}\begin{array}{c}\IA\\\downarrow\\\IA\end{array},\varphi_{\IA}\right)$ is a relative Rota-Baxter operator where $$\varphi_{\IA}=p_{A}\circ\varphi_{A}^{ad}\circ (i_{A}\otimes i_{A}).$$
\begin{figure}[h]
\[
\xymatrix{&\IA\ar[rr]^-{\lambda_{\IA}}\ar[d]^-{i_{A}} & &\IA\ar[d]^-{i_{A}} \\&A\ar@<-1ex>[d]_-{g}\ar[rr]^-{\lambda_{A}} & &A\ar@<-1ex>[d]_-{g}\\&H\ar@(dr,dl)[]\ar@<-1ex>[u]_-{f}\ar[rr]^-{\lambda_{H}}& &H.\ar@(dr,dl)[]\ar@<-1ex>[u]_-{f}}
\]
\caption{$\lambda_{\IA}$ is a relative Rota-Baxter operator.}
\end{figure}
\end{example}

\begin{corollary}
There exists a functor 
\[{\sf P}\colon {\sf P({\sf coc}\textnormal{-}{\sf rRB})}\longrightarrow {\sf coc}\textnormal{-}{\sf rRB}\]
acting on objects by 
\[{\sf P}\left(\left(\left(T\begin{array}{c}H\\\downarrow\\B \end{array},\varphi_{H}\right),\left(L\begin{array}{c}A\\\downarrow\\D \end{array},\varphi_{A}\right),f,h,g,l \right)\right)=\left(L_{0}\begin{array}{c}\IA\\\downarrow\\\ID \end{array},\varphi_{\IA}\right)\]
and on morphisms by ${\sf P}((x,y,z,t))=(y_{0},t_{0})$, where $y_{0}\colon\IA\rightarrow\IAp$ is the unique morphism satisfying that $i_{A'}\circ y_{0}=y\circ i_{A}$ and $t_{0}\colon\ID\rightarrow\IDp$ is the unique morphism such that $i_{D'}\circ t_{0}=t\circ i_{D}$.
\end{corollary}
\begin{proof}
Functor ${\sf P}$ is well-defined on objects thanks to Theorem \ref{thT0}. Consider $(x,y,z,t)$ a morphism in ${\sf P({\sf coc}\textnormal{-}{\sf rRB})}$ between $$\left(\left(T\begin{array}{c}H\\\downarrow\\B \end{array},\varphi_{H}\right),\left(L\begin{array}{c}A\\\downarrow\\D \end{array},\varphi_{A}\right),f,h,g,l\right)\;\; {\rm and}\;\;\left(\left(T'\begin{array}{c}H'\\\downarrow\\B' \end{array},\varphi_{H'}\right),\left(L'\begin{array}{c}A'\\\downarrow\\D' \end{array},\varphi_{A'}\right),f',h',g',l' \right).$$
On the one hand, note that $y\circ i_{A}$ factors through the equalizer $i_{A'}$. Indeed,
\begin{align*}
&(A'\otimes g')\circ\delta_{A'}\circ y\circ i_{A}\\=&(y\otimes(g'\circ y))\circ\delta_{A}\circ i_{A}\;\footnotesize\textnormal{(by the condition of coalgebra morphism for $y$)}\\=&(y\otimes(x\circ g))\circ\delta_{A}\circ i_{A}\;\footnotesize\textnormal{(by \eqref{proj-morph} for $(x,y)$)}\\=&(y\circ i_{A})\otimes(x\circ\eta_{H})\;\footnotesize\textnormal{(by the equalizer condition for $i_{A}$)}\\=&(y\circ i_{A})\otimes\eta_{H'}\;\footnotesize\textnormal{(by the condition of algebra morphism for $x$)}.
\end{align*}

Therefore, there exists an unique $y_{0}\colon \IA\rightarrow \IAp$ satisfying that
\begin{equation}\label{xBcons}
i_{A'}\circ y_{0}=y\circ i_{A}.
\end{equation}
As a consequence, if we compose on the left with $p_{A'}$, then the following equality is obtained:
\begin{equation}\label{xBdef}
y_{0}=p_{A'}\circ y\circ i_{A}.
\end{equation} 
Note also that $y_{0}$ is a Hopf algebra morphism. On the one side, it is an algebra morphism because
\begin{align*}
&i_{A'}\circ y_{0}\circ\eta_{\IA}\\=&y\circ i_{A}\circ \eta_{\IA}\;\footnotesize\textnormal{(by \eqref{xBcons})}\\=&y\circ \eta_{A}\;\footnotesize\textnormal{(by the equality $i_{A}\circ\eta_{\IA}=\eta_{A}$)}\\=&\eta_{A'}\;\footnotesize\textnormal{(by the condition of algebra morphism for $y$)}\\=&i_{A'}\circ\eta_{\IAp}\;\footnotesize\textnormal{(by the equality $i_{A'}\circ\eta_{\IAp}=\eta_{A'}$)}
\end{align*}
and also
\begin{align*}
&i_{A'}\circ y_{0}\circ \mu_{\IA}\\=&y\circ i_{A}\circ \mu_{\IA}\;\footnotesize\textnormal{(by \eqref{xBcons})}\\=&y\circ\mu_{A}\circ(i_{A}\otimes i_{A})\;\footnotesize\textnormal{(by \eqref{m-i})}\\=&\mu_{A'}\circ((y\circ i_{A})\otimes(y\circ i_{A}))\;\footnotesize\textnormal{(by the condition of algebra morphism for $y$)}\\=&\mu_{A'}\circ(i_{A'}\otimes i_{A'})\circ(y_{0}\otimes y_{0})\;\footnotesize\textnormal{(by \eqref{xBcons})}\\=&i_{A'}\circ\mu_{\IAp}\circ(y_{0}\otimes y_{0})\;\footnotesize\textnormal{(by \eqref{m-i})}.
\end{align*}

On the other side, $\varepsilon_{\IAp}\circ y_{0}=\varepsilon_{\IA}$ follows by \eqref{xBcons} and the fact that $y$ preserves the counit and
\begin{align*}
&\delta_{\IAp}\circ y_{0}\\=&(p_{A'}\otimes p_{A'})\circ\delta_{A'}\circ y\circ i_{A}\;\footnotesize\textnormal{(by \eqref{xBcons})}\\=&(p_{A'}\otimes p_{A'})\circ(y\otimes y)\circ \delta_{A}\circ i_{A}\;\footnotesize\textnormal{(by the condition of coalgebra morphism for $y$)}\\=&(y_{0}\otimes y_{0})\circ\delta_{\IA}\;\footnotesize\textnormal{(by the condition of coalgebra morphism for $i_{A}$ and \eqref{xBdef})},
\end{align*}
so $y_{0}$ is a coalgebra morphism too. Analogously, it results easy to prove that $t\circ i_{D}$ factors through the equalizer $i_{D'}$, and then there exists an unique $t_{0}\colon \ID\rightarrow \IDp$ such that the equalities
\begin{gather}\label{zCcons}
i_{D'}\circ t_{0}=t\circ i_{D},\\\label{zCdef}t_{0}=p_{D'}\circ t\circ i_{D},
\end{gather}
hold and $t_{0}$ is a Hopf algebra morphism.

Thus, to conclude the proof it is enough to see that the pair $(y_{0},t_{0})$ satisfies conditions \eqref{cond1morrRB} and \eqref{cond2morrRB}. At first, \eqref{cond1morrRB} follows by
\begin{align*}
&i_{D'}\circ L_{0}'\circ y_{0}\\=&L'\circ i_{A'}\circ y_{0}\;\footnotesize\textnormal{(by \eqref{T0cons})}\\=&L'\circ y\circ i_{A}\;\footnotesize\textnormal{(by \eqref{xBcons})}\\=&t\circ L\circ i_{A}\;\footnotesize\textnormal{(by \eqref{cond1morrRB} for $(y,t)$)}\\=&t\circ i_{D}\circ L_{0}\;\footnotesize\textnormal{(by \eqref{T0cons})}\\=&i_{D'}\circ t_{0}\circ L_{0}\;\footnotesize\textnormal{(by \eqref{zCcons})}.
\end{align*}

In addition, \eqref{cond2morrRB} is consequence of 
\begin{align*}
&i_{A'}\circ y_{0}\circ \varphi_{\IA}\\=&y\circ i_{A}\circ\varphi_{\IA}\;\footnotesize\textnormal{(by \eqref{xBcons})}\\=&y\circ \varphi_{A}\circ (i_{D}\otimes i_{A})\;\footnotesize\textnormal{(by \eqref{PhiIBcons})}\\=&\varphi_{A'}\circ ((t\circ i_{D})\otimes(y\circ i_{A}))\;\footnotesize\textnormal{(by \eqref{cond2morrRB} for $(y,t)$)}\\=&\varphi_{A'}\circ(i_{D'}\otimes i_{A'})\circ(t_{0}\otimes y_{0})\;\footnotesize\textnormal{(by \eqref{xBcons} and \eqref{zCcons})}\\=&i_{A'}\circ \varphi_{\IAp}\circ(t_{0}\otimes y_{0})\;\footnotesize\textnormal{(by \eqref{PhiIBcons})}.
\end{align*}

Then, ${\sf P}$ is also well-defined on morphisms.

\begin{figure}[h]
\begin{minipage}{0.5 \linewidth}
\[\hspace{-2.5 cm}\xymatrix{&\IA\ar@{-->}[d]^-{\exists^{\bullet}\,y_{0}}\ar[rr]^-{i_{A}} & &A\ar[d]^-{y}\ar@<1ex>[rr]^-{(A\otimes g)\circ\delta_{A}}\ar@<-1ex>[rr]_-{A\otimes\eta_{H}} & &A\otimes H\\&\IAp \ar[rr]^-{i_{A'}} & &A'\ar@<1ex>[rr]^-{(A'\otimes g')\circ\delta_{A'}}\ar@<-1ex>[rr]_-{A' \otimes\eta_{H'}} & &A'\otimes H'}\]
\end{minipage}
\begin{minipage}{0.4 \linewidth}
\[\hspace{-1.5 cm}\xymatrix{&\ID\ar@{-->}[d]^-{\exists^{\bullet}\,t_{0}}\ar[rr]^-{i_{D}} & &D\ar[d]^-{t}\ar@<1ex>[rr]^-{(D\otimes l)\circ\delta_{D}}\ar@<-1ex>[rr]_-{D\otimes\eta_{B}} & &D\otimes B\\&\IDp \ar[rr]^-{i_{D'}} & &D'\ar@<1ex>[rr]^-{(D'\otimes l')\circ\delta_{D'}}\ar@<-1ex>[rr]_-{D' \otimes\eta_{B'}} & &D'\otimes B'}\]
\end{minipage}
\caption{The construction of $y_{0}$ and $t_{0}$.}
\end{figure}
\end{proof}

In what follows we recall the notion of projection of Hopf braces introduced in \cite{FGRR}.

\begin{definition}
A projection of Hopf braces in {\sf C} is a 4-tuple  $({\mathbb H}, \mathbb{D}, x,y)$, where ${\mathbb H}$, $\mathbb{D} $  are Hopf braces in {\sf C}, $x:{\mathbb H}\rightarrow \mathbb{D}$, $y:\mathbb{D}\rightarrow {\mathbb H}$ are morphisms of Hopf braces in ${\sf C}$ and the following equality $y\circ x=id_{\mathbb H}$ holds.

A morphism between two projections of Hopf braces  $({\mathbb H}, \mathbb{D}, x,y)$ and  $({\mathbb H}^{\prime}, \mathbb{D}^{\prime}, x^{\prime},y^{\prime})$ is a pair $(z ,t)$ where $z:{\mathbb H}\rightarrow {\mathbb H}^{\prime}$,  $t:\mathbb{D}\rightarrow \mathbb{D}^{\prime}$ are morphisms in {\sf HBr} and the following equalities hold:
\begin{equation}
\label{eqpHBr}
x^{\prime}\circ z=t\circ x, \;\;\;\; y^{\prime}\circ t=z\circ y.
\end{equation}
	
With these morphisms and the previous objects we can define the category of projections of Hopf braces. We will denote this category by {\sf P(HBr)}. With ${\sf P}({\sf coc}\textnormal{-}{\sf  HBr})$ we will denote the category of projections between cocommutative Hopf braces.

\end{definition}
	
\begin{remark} Following \cite{FGRR}, if $({\mathbb H}, {\mathbb D}, x,y)$ is a projection of Hopf braces in {\sf C}, then we have two projections of Hopf algebras $(H_{1}, D_{1}, x,y)$ and $(H_{2}, D_{2}, x,y)$. Then, with $q_{D}^1$ and $q_{D}^2$ we will denote the associated idempotent morphisms.  Note that, if $q_{D}^1=i_{D}^1\circ p_{D}^1 $ and $q_{D}^2=i_{D}^2\circ p_{D}^2$, with  $ p_{D}^1 \circ  i_{D}^1 =id_{I(q_{D}^1)}$ and $p_{D}^2 \circ  i_{D}^2 =id_{I(q_{D}^2)}$, we have that 
$$
\xymatrix{&I(q_{D}^{k})\ar[rr]^-{i_{D}^{k}} & &D\ar@<1ex>[rr]^-{(D\otimes y)\circ \delta_{D}}\ar@<-1ex>[rr]_-{D\otimes\eta_{H}} & &D\otimes H}
$$
is an equalizer diagram for $k\in \{1,2\} $ and, as a consequence, we can assume that  $i_{D}^1=i_{D}^2=i_{D}$  and $I(q_{D}^1)=I(q_{D}^2)=\ID$.  Then,  $p_{D}^1\circ i_{D}=id_{I(q_{D})}=p_{D}^2\circ i_{D}$ holds.
\end{remark}

Just as there exists an adjunction between the category of relative Rota-Baxter operators, ${\sf rRB}^{\star}$, and the category of cocommutative Hopf braces, {\sf coc-HBr} (see Theorem \ref{adjunction}), it seems reasonable to ask whether this adjunction carries over to their respective projection categories.

\begin{theorem}
There exists a functor 
\[{\sf Q}\colon {\sf P(coc\textnormal{-}HBr)}\longrightarrow {\sf P(rRB^{\star})}\]
acting on objects by ${\sf Q}((\mathbb{H},\mathbb{D},x,y))=({\sf F}(\mathbb{H}),{\sf F}(\mathbb{B}),x,x,y,y)$ and on morphisms by ${\sf Q}((z,t))=(z,t,z,t)$, where ${\sf F}$ denotes the functor introduced in Theorem \ref{adjunction}.
\end{theorem}
\begin{proof} It follows by the fact that functor ${\sf F}\colon {\sf coc\textnormal{-}HBr}\longrightarrow {\sf rRB^{\star}}$ introduced in Theorem \ref{adjunction} is well-defined.
\end{proof}

\begin{remark}
Note that the image of ${\sf Q}$ is in 
${\sf P(coc\textnormal{-}rRB)}$.  Then if 
$I^{\star}:{\sf P({\sf coc\textnormal{-}rRB})}\longrightarrow {\sf P({\sf rRB^{\star}})}$ denotes the inclusion functor, we have that ${\sf Q}=I^{\star}\circ 
{\sf Q}'$ where \[{\sf Q}'\colon {\sf P(coc\textnormal{-}HBr)}\longrightarrow {\sf P(coc\textnormal{-}rRB)}\] is the functor defined as ${\sf Q}$.
\end{remark}

\begin{theorem}
\label{RRR}
There exists a functor 
\[{\sf R}\colon {\sf P(rRB^{\star})}\longrightarrow {\sf P(coc\textnormal{-}HBr)}\]
acting on objects by 
\begin{align*} &{\sf R}\left(\left(\left(T\begin{array}{c}H\\\downarrow\\B \end{array},\varphi_{H}\right),\left(L\begin{array}{c}A\\\downarrow\\D \end{array},\varphi_{A}\right),f,h,g,l\right)\right)\\=&\left({\sf G}\left(\left(T\begin{array}{c}H\\\downarrow\\B \end{array},\varphi_{H}\right)\right),{\sf G}\left(\left(L\begin{array}{c}A\\\downarrow\\D\end{array},\varphi_{A}\right)\right),f,g\right), \end{align*}
where ${\sf G}$ denotes the functor introduced in Theorem \ref{adjunction}, and on morphisms by ${\sf R}((x,y,z,t))=(x,y)$.
\end{theorem}

\begin{proof}
The proof is straightforward taking into account that functor ${\sf G}\colon {\sf rRB^{\star}}\longrightarrow {\sf coc\textnormal{-}HBr}$ is well-defined, as it is possible to consult in Theorem \ref{adjunction}.
\end{proof}

\begin{theorem}\label{projadjuncQR}
The functor ${\sf Q}$ is left adjoint of functor ${\sf R}$.
\end{theorem}
\begin{proof}
We have to show that there exists a bijection ${}^{(\mathbb{H},\mathbb{D},x,y)}\Sigma_{(T,L,f,g,h,l)}$ between

{\small 
\begin{align*}
&\Hom_{{\sf P(rRB^{\star})}}\left({\sf Q}\left(\left(\mathbb{H},\mathbb{D},x,y\right)\right),\left(\left(T\begin{array}{c}X\\\downarrow\\B \end{array},\varphi_{X}\right),\left(L\begin{array}{c}A\\\downarrow\\Y \end{array},\varphi_{A}\right),f,h,g,l\right)\right)\\=&\Hom_{{\sf P(rRB^{\star})}}\left(\left(\left(id_{H}\begin{array}{c}H_{1}\\\downarrow\\H_{2}\end{array},\Gamma_{H_{1}}\right),\left(id_{D}\begin{array}{c}D_{1}\\\downarrow\\D_{2}\end{array},\Gamma_{D_{1}}\right),x,x,y,y\right), \left(\left(T\begin{array}{c}X\\\downarrow\\B \end{array},\varphi_{X}\right),\left(L\begin{array}{c}A\\\downarrow\\Y \end{array},\varphi_{A}\right),f,h,g,l\right)\right)
\end{align*}
}
and
\small{
\begin{align*}
&\Hom_{{\sf P(coc\textnormal{-}HBr)}}\left((\mathbb{H},\mathbb{D},x,y),{\sf R}\left(\left(\left(T\begin{array}{c}X\\\downarrow\\B \end{array},\varphi_{X}\right),\left(L\begin{array}{c}A\\\downarrow\\Y \end{array},\varphi_{A}\right),f,h,g,l\right)\right)\right)\\=&
\Hom_{{\sf P(coc\textnormal{-}HBr)}}\left((\mathbb{H},\mathbb{D},x,y),(\overline{\mathbb{X}},\overline{\mathbb{A}},f,g)\right).
\end{align*}}
	
At first, consider $(a,b,c,d)$ a morphism of relative Rota-Baxter projections between
\[\left(\left(id_{H}\begin{array}{c}H_{1}\\\downarrow\\H_{2}\end{array},\Gamma_{H_{1}}\right),\left(id_{D}\begin{array}{c}D_{1}\\\downarrow\\D_{2}\end{array},\Gamma_{D_{1}}\right),x,x,y,y\right)\;{\rm and}\;\left(\left(L\begin{array}{c}X\\\downarrow\\B \end{array},\varphi_{X}\right),\left(L\begin{array}{c}A\\\downarrow\\Y \end{array},\varphi_{Y}\right),f,h,g,l\right),\]
and let's show that $(a,b)$ is a morphism in ${\sf P(coc\textnormal{-}HBr)}$ between $(\mathbb{H},\mathbb{D},x,y)$ and $(\overline{\mathbb{X}},\overline{\mathbb{A}},f,g)$. Indeed, we only need to show that $a$, $b$ are multiplicative morphisms. Then, on the one hand
\begin{align*}
&\overline{\mu}_{X}\circ(a\otimes a)\\=&\mu_{X}\circ(a\otimes(\varphi_{X}\circ ((T\circ a)\otimes a)))\circ(\delta_{H}\otimes H)\;\footnotesize\textnormal{(by the condition of coalgebra morphism for $a$)}\\=&\mu_{X}\circ(a\otimes(\varphi_{X}\circ (c\otimes a)))\circ(\delta_{H}\otimes H)\;\footnotesize\textnormal{(by \eqref{cond1morrRB} for $(a,c)$)}\\=&\mu_{X}\circ(a\otimes (a\circ\Gamma_{H_{1}}))\circ(\delta_{H}\otimes H)\;\footnotesize\textnormal{(by \eqref{cond2morrRB} for $(a,c)$)}\\=&a\circ\mu_{H}^{1}\circ(H\otimes\Gamma_{H_{1}})\circ(\delta_{H}\otimes H)\;\footnotesize\textnormal{(by the condition of algebra morphism for $a\colon H_{1}\rightarrow X$)}\\=&a\circ\mu_{H}^{2}\;\footnotesize\textnormal{(by \eqref{eb2})},
\end{align*}
and, on the other hand, $\overline{\mu}_{A}\circ(b\otimes b)=b\circ\mu_{D}^{2}$ following the same arguments as before. Then, we define
\[{}^{(\mathbb{H},\mathbb{D},x,y)}\Sigma_{(T,L,f,g,h,l)}(a,b,c,d)=(a,b).\]

Now, let $(z,t)$ be a morphism of projections of Hopf braces between $(\mathbb{H},\mathbb{D},x,y)$ and $(\overline{\mathbb{X}},\overline{\mathbb{A}},f,g)$. Then, let's prove that $(z,t,T\circ z,L\circ t)$ is a morphism of projections of relative Rota-Baxter operators between 
\[\left(\left(id_{H}\begin{array}{c}H_{1}\\\downarrow\\H_{2}\end{array},\Gamma_{H_{1}}\right),\left(id_{D}\begin{array}{c}D_{1}\\\downarrow\\D_{2}\end{array},\Gamma_{D_{1}}\right),x,x,y,y\right)\;{\rm and}\;\left(\left(T\begin{array}{c}X\\\downarrow\\B \end{array},\varphi_{X}\right),\left(L\begin{array}{c}A\\\downarrow\\Y \end{array},\varphi_{A}\right),f,h,g,l\right).\]
	
By hypothesis, $(z,t)\colon (H_{1},D_{1},x,y)\rightarrow (X,A,f,g)$ is a morphism of projections of Hopf algebras. On the other hand, $(T\circ z, L\circ t)\colon (H_{2},D_{2},x,y)\rightarrow (B,Y,h,l)$ is a morphism of projections of Hopf algebras too. That is due to the fact that $T\circ z$ and $L\circ t$ are Hopf algebra morphisms, what follows by
\begin{align*}
&\mu_{B}\circ((T\circ z)\otimes(T\circ z))\\=&T\circ\overline{\mu}_{X}\circ(z\otimes z)\;\footnotesize\textnormal{(by \eqref{rRBcond} for $T$)}\\=&T\circ z\circ\mu_{H}^{2}\;\footnotesize\textnormal{(by the condition of algebra morphism for $z\colon H_{2}\rightarrow \overline{X}$)}
\end{align*}
and
\begin{align*}
&\mu_{Y}\circ((L\circ t)\otimes(L\circ t))\\=&L\circ\overline{\mu}_{A}\circ(t\otimes t)\;\footnotesize\textnormal{(by \eqref{rRBcond} for $L$)}\\=&L\circ t\circ \mu_{D}^{2}\;\footnotesize\textnormal{(by the condition of algebra morphism for $t\colon D_{2}\rightarrow \overline{A}$)},
\end{align*}
and also we have that
\begin{align*}
&h\circ T\circ z\\=&L\circ f\circ z\;\footnotesize\textnormal{(by \eqref{cond1morrRB} for $(f,h)$)}\\=&L\circ t\circ x\;\footnotesize\textnormal{(by \eqref{proj-morph} for $(z,t)$)},
\end{align*}
and
\begin{align*}
&l\circ L\circ t\\=&T\circ g\circ t\;\footnotesize\textnormal{(by \eqref{cond1morrRB} for $(g,l)$)}\\=&T\circ z\circ y\;\footnotesize\textnormal{(by \eqref{proj-morph} for $(z,t)$)}.
\end{align*}

Note also that $(z,T\circ z)$ is a morphism of relative Rota-Baxter operators between $\left(id_{H}\begin{array}{c}H_{1}\\\downarrow\\H_{2}\end{array},\Gamma_{H_{1}}\right)$ and $\left(T\begin{array}{c}X\\\downarrow\\B\end{array},\varphi_{X}\right)$ because \eqref{cond1morrRB} follows directly and 
\begin{align*}
&z\circ\Gamma_{H_{1}}\\=&\mu_{X}\circ ((z\circ\lambda_{H}^{1})\otimes(z\circ\mu_{H}^{2}))\circ(\delta_{H}\otimes H)\;\footnotesize\textnormal{(by the condition of algebra morphism for $z\colon H_{1}\rightarrow X$)}\\=&\mu_{X}\circ(\lambda_{X}\otimes\overline{\mu}_{X})\circ(((z\otimes z)\circ\delta_{H})\otimes z)\;\footnotesize\textnormal{(by \eqref{morant} and the condition of algebra morphism for $z\colon H_{2}\rightarrow \overline{X}$)}\\=&\mu_{X}\circ ((\lambda_{X}\ast id_{X})\otimes(\varphi_{X}\circ(T\otimes X)))\circ (\delta_{X}\otimes X)\circ(z\otimes z)\;\footnotesize\textnormal{(by the condition of coalgebra morphism for $z$,}\\&\footnotesize\textnormal{the associativity of $\mu_{X}$ and coassociativity of $\delta_{X}$)}\\=&\varphi_{X}\circ((T\circ z)\otimes z)\;\footnotesize\textnormal{(by \eqref{antipode} and the (co)unit property)}.
\end{align*}

Following the same arguments, we can also show that $(t,L\circ t)$ is a morphism of relative Rota-Baxter operators between $\left(id_{D}\begin{array}{c}D_{1}\\\downarrow\\D_{2}\end{array},\Gamma_{D_{1}}\right)$ and $\left(L\begin{array}{c}A\\\downarrow\\Y\end{array},\varphi_{A}\right)$, so we define
\[({}^{(\mathbb{H},\mathbb{D},x,y)}\Sigma_{(T,L,f,h,g,l)})^{-1}(z,t)=(z,t,T\circ z,L\circ t).\]
	
Therefore, it only remains to compute that ${}^{(\mathbb{H},\mathbb{D},x,y)}\Sigma_{(T,L,f,g,h,l)}$ and $({}^{(\mathbb{H},\mathbb{D},x,y)}\Sigma_{(T,L,f,g,h,l)})^{-1}$ define a bijection. Indeed,
\begin{align*}
&(({}^{(\mathbb{H},\mathbb{D},x,y)}\Sigma_{(T,L,f,g,h,l)})^{-1}\circ{}^{(\mathbb{H},\mathbb{D},x,y)}\Sigma_{(T,L,f,g,h,l)})(a,b,c,d)\\=&({}^{(\mathbb{H},\mathbb{D},x,y)}\Sigma_{(T,L,f,g,h,l)})^{-1}(a,b)=(a,b,T\circ a, L\circ b)=(a,b,c,d)
\end{align*}
because $T\circ a=c$ and $L\circ b=d$ by \eqref{cond1morrRB} for $(a,c)$ and $(b,d)$, and also
\begin{align*}
&({}^{(\mathbb{H},\mathbb{D},x,y)}\Sigma_{(T,L,f,g,h,l)}\circ({}^{(\mathbb{H},\mathbb{D},x,y)}\Sigma_{(T,L,f,g,h,l)})^{-1})(z,t)={}^{(\mathbb{H},\mathbb{D},x,y)}\Sigma_{(T,L,f,g,h,l)}(z,t,T\circ z,L\circ t)=(z,t).\qedhere
\end{align*}
\end{proof}

The following result, proved in \cite{FGRR}, asserts that  every projection of cocommutative Hopf braces give rise to a new Hopf brace in {\sf C}.

\begin{theorem}
If $(\mathbb{H},\mathbb{D},x,y)\in{\sf P({\sf coc\textnormal{-}HBr})}$, then 
\[\mathbb{I}(q_{D})=(\ID,\eta_{\ID},\mu_{\ID}^{1},\mu_{\ID}^{2},\varepsilon_{\ID},\delta_{\ID},\lambda_{\ID}^{1},\lambda_{\ID}^{2})\]
is a Hopf brace in {\sf C}, where
\begin{gather}\label{etaIB}
\eta_{\ID}=p_{D}^{1}\circ\eta_{D}=p_{D}^{2}\circ\eta_{D},\\\label{mu1IB}
\mu_{\ID}^{1}=p_{D}^{1}\circ\mu_{D}^{1}\circ (i_{D}\otimes i_{D})=p_{D}^{2}\circ\mu_{D}^{1}\circ(i_{D}\otimes i_{D}),\\\label{mu2IB}
\mu_{\ID}^{2}=p_{D}^{2}\circ\mu_{D}^{2}\circ(i_{D}\otimes i_{D})=p_{D}^{1}\circ\mu_{D}^{2}\circ(i_{D}\otimes i_{D}),\\\label{epsIB}
\varepsilon_{\ID}=\varepsilon_{D}\circ i_{D},\\\label{deltaIB}
\delta_{\ID}=(p_{D}^{1}\otimes p_{D}^{1})\circ\delta_{D}\circ i_{D}=(p_{D}^{2}\otimes p_{D}^{2})\circ\delta_{D}\circ i_{D},\\\label{lambda1IB}
\lambda_{\ID}^{1}=p_{D}^{1}\circ\lambda_{D}^{1}\circ i_{D},\\\label{lambda2IB}
\lambda_{\ID}^{2}=p_{D}^{2}\circ\lambda_{D}^{2}\circ i_{D}.
\end{gather}
\end{theorem}

\begin{corollary}
There exists a functor 
\[{\sf P}'\colon {\sf P({\sf coc\textnormal{-}HBr})}\longrightarrow {\sf coc\textnormal{-}HBr}\]
acting on objects by ${\sf P}'((\mathbb{H},\mathbb{D},x,y))=\mathbb{I}(q_{D})$ and on morphism by ${\sf P}'((z,t))=t_{0}$, where $t_{0}\colon \ID\rightarrow I(q_{D'})$ is the unique morphism verifying that $i_{D'}\circ t_{0}=t\circ i_{D}$.
\end{corollary}
\begin{proof}
Thanks to previous theorem, functor ${\sf P}'$ is well-defined on objects. Let's show that it is well-defined on morphisms too. Take $(z,t)\colon (\mathbb{H},\mathbb{D},x,y)\rightarrow (\mathbb{H}',\mathbb{D}',x',y')$ a morphism of projections of Hopf braces and note that $t\circ i_{D}$ factors through the equalizer $i_{D'}$:
\begin{align*}
&(D'\otimes y')\circ\delta_{D'}\circ t\circ i_{D}\\=&(t\otimes (y'\circ t))\circ\delta_{D}\circ i_{D}\;\footnotesize\textnormal{(by the condition of coalgebra morphism for $t$)}\\=&(t\otimes (z\circ y))\circ\delta_{D}\circ i_{D}\;\footnotesize\textnormal{(by \eqref{eqpHBr})}\\=&(t\circ i_{D})\otimes(z\circ \eta_{H})\;\footnotesize\textnormal{(by the equalizer condition for $i_{D}$)}\\=&(t\circ i_{D})\otimes\eta_{H'}\;\footnotesize\textnormal{(by the condition of algebra morphism for $z$)}.
\end{align*}

Therefore, there exists an unique $t_{0}\colon \ID\rightarrow \IDp$ such that
\begin{equation}\label{tBcons}
i_{D'}\circ t_{0}=t\circ i_{D},
\end{equation}
what implies that
\begin{equation}\label{tBdef}
t_{0}=p_{D'}^{1}\circ t\circ i_{D}=p_{D'}^{2}\circ t\circ i_{D}.
\end{equation}
\begin{figure}[h]
\[\xymatrix{&\ID\ar@{-->}[d]^-{\exists^{\bullet}\,t_{0}}\ar[rr]^-{i_{D}} & &D\ar[d]^-{t}\ar@<1ex>[rr]^-{(D\otimes y)\circ\delta_{D}}\ar@<-1ex>[rr]_-{D\otimes\eta_{H}} & &D\otimes H\\&\IDp \ar[rr]^-{i_{D'}} & &D'\ar@<1ex>[rr]^-{(D'\otimes y')\circ\delta_{D'}}\ar@<-1ex>[rr]_-{D' \otimes\eta_{H'}} & &D'\otimes H'}\]
\label{tBconsfigure}
\caption{The construction of $t_{0}$.}
\end{figure}

Then, to conclude it is enough to show that $t_{0}\colon\mathbb{I}(q_{D})\rightarrow\mathbb{I}(q_{D'})$ is a morphism of Hopf braces. The condition of coalgebra morphism for $t_{0}$ follows by
\begin{align*}
&\delta_{\IDp}\circ t_{0}\\=&(p_{D'}^{1}\otimes p_{D'}^{1})\circ\delta_{D'}\circ t\circ i_{D}\;\footnotesize\textnormal{(by \eqref{tBcons})}\\=&((p_{D'}^{1}\circ t)\otimes (p_{D'}^{1}\circ t))\circ\delta_{D}\circ i_{D}\;\footnotesize\textnormal{(by the condition of coalgebra morphism for $t$)}\\=&(t_{0}\otimes t_{0})\circ\delta_{\ID}\;\footnotesize\textnormal{(by the condition of coalgebra morphism for $i_{D}$ and \eqref{tBdef})},
\end{align*}
and also $t_{0}$ is an algebra morphism between $I(q_{D})_{k}$ and $I(q_{D'})_{k}$ for all $k=1,2$ because
\begin{align*}
&i_{D'}\circ t_{0}\circ\mu_{\ID}^{k}\\=&t\circ i_{D}\circ \mu_{\ID}^{k}\;\footnotesize\textnormal{(by \eqref{tBcons})}\\=&t\circ\mu_{D}^{k}\circ (i_{D}\otimes i_{D})\;\footnotesize\textnormal{(by \eqref{m-i})}\\=&\mu_{D'}^{k}\circ((t\circ i_{D})\otimes(t\circ i_{D}))\;\footnotesize\textnormal{(by the condition of morphism of Hopf braces for $t$)}\\=&\mu_{D'}^{k}\circ(i_{D'}\otimes i_{D'})\circ(t_{0}\otimes t_{0})\;\footnotesize\textnormal{(by \eqref{tBcons})}\\=&i_{D'}\circ\mu_{\IDp}^{k}\circ(t_{0}\otimes t_{0})\;\footnotesize\textnormal{(by \eqref{m-i})}.\qedhere
\end{align*}
\end{proof}

\begin{theorem}\label{commutativeDiag1}
The diagram of functors 
\[\xymatrix{
&{\sf P(coc\textnormal{-}rRB)}\ar[d]^-{{\sf P}}\ar[rr]^-{{\sf R}'} & &{\sf P(coc\textnormal{-}HBr)}\ar[d]^-{{\sf P}'}\\
&{\sf coc\textnormal{-}rRB}\ar[rr]^-{{\sf G}} & &{\sf coc\textnormal{-}HBr}
}\]
where ${\sf R}'$ is the restriction of the functor ${\sf R}\colon {\sf P(rRB^{\star})}\longrightarrow {\sf P(coc\textnormal{-}HBr)}$, introduced in Theorem \ref{RRR}, to the subcategory ${\sf P(coc\textnormal{-}rRB)}$, is commutative. 
\end{theorem}
\begin{proof}

We only detail that the previous diagram is commutative on objects because the commutativity on morphisms is straightforward. Let $\left(\left(T\begin{array}{c}H\\\downarrow\\B \end{array},\varphi_{H}\right),\left(L\begin{array}{c}A\\\downarrow\\D \end{array},\varphi_{A}\right),f,h,g,l\right)$ be a projection of cocommutative relative Rota-Baxter operators. Therefore, we have to show that the Hopf brace
\[{\sf G}\left(\left(L_{0}\begin{array}{c}\IA\\\downarrow\\\ID\end{array},\varphi_{\IA}\right)\right)=\overline{\mathbb{I}(q_{A})}=(I(q_{A}),\eta_{\IA},\mu_{\IA},\overline{\mu}_{\IA},\varepsilon_{\IA},\delta_{\IA},\lambda_{\IA},\overline{\lambda}_{\IA}),\]
where $\eta_{\IA}$, $\mu_{\IA}$, $\varepsilon_{\IA}$, $\delta_{\IA}$ and $\lambda_{\IA}$ are the usual ones, while 
\begin{gather}
\overline{\mu}_{\IA}=\mu_{\IA}\circ(\IA\otimes (\varphi_{\IA}\circ(L_{0}\otimes\IA)))\circ(\delta_{\IA}\otimes\IA),\\
\overline{\lambda}_{\IA}=\varphi_{\IA}\circ((\lambda_{\ID}\circ L_{0})\otimes \lambda_{\IA})\circ\delta_{\IA},
\end{gather}
is the same as the Hopf brace 
\[{\sf P}'(\overline{\mathbb{H}},\overline{\mathbb{A}},f,g)=\mathbb{I}(q_{\overline{A}})=(\IA,\eta_{I(q_{\overline{A}})},\mu_{I(q_{\overline{A}})}^{1},\mu_{I(q_{\overline{A}})}^{2},\varepsilon_{I(q_{\overline{A}})},\delta_{I(q_{\overline{A}})},\lambda_{I(q_{\overline{A}})}^{1},\lambda_{I(q_{\overline{A}})}^{2}),\]
where $\eta_{I(q_{\overline{A}})}=\eta_{\IA}$, $\mu_{I(q_{\overline{A}})}^{1}=\mu_{\IA}$, $\varepsilon_{I(q_{\overline{A}})}=\varepsilon_{\IA}$, $\delta_{I(q_{\overline{A}})}=\delta_{\IA}$ and $\lambda_{I(q_{\overline{A}})}^{1}=\lambda_{I(q_{\overline{A}})}$, while $\mu_{I(q_{\overline{A}})}^{2}$ is the unique product satisfying that 
\begin{equation}\label{overB2mu}i_{A}\circ \mu_{I(q_{\overline{A}})}^{2}=\overline{\mu}_{A}\circ (i_{A}\otimes i_{A})\end{equation} and $\lambda_{I(q_{\overline{A}})}^{2}$ is the unique morphism which verify that 
\begin{equation*}i_{A}\circ \lambda_{I(q_{\overline{A}})}^{2}=\overline{\lambda}_{A}\circ i_{A}.\end{equation*} 

Thus, to conclude the proof it is enough to show that $\overline{\mu}_{\IA}=\mu^{2}_{I(q_{\overline{A}})}$ what implies that $\overline{\lambda}_{\IA}=\lambda_{I(q_{\overline{A}})}^{2}$ due the uniqueness of the antipode for a bialgebra structure. Indeed,
\begin{align*}
&i_{A}\circ \overline{\mu}_{\IA}\\=&\mu_{A}\circ (i_{A}\otimes(i_{A}\circ\varphi_{\IA}\circ (L_{0}\otimes \IA)))\circ(\delta_{\IA}\otimes\IA)\;\footnotesize\textnormal{(by \eqref{m-i})}\\=&\mu_{A}\circ(i_{A}\otimes (\varphi_{A}\circ ((i_{D}\circ L_{0})\otimes i_{A})))\circ(\delta_{\IA}\otimes\IA)\;\footnotesize\textnormal{(by \eqref{PhiIBcons})}\\=&\mu_{A}\circ(i_{A}\otimes (\varphi_{A}\circ (L\otimes A)))\circ(((i_{A}\otimes i_{A})\circ\delta_{\IA})\otimes i_{A})\;\footnotesize\textnormal{(by \eqref{T0cons})}\\=&\overline{\mu}_{A}\circ(i_{A}\otimes i_{A})\;\footnotesize\textnormal{(by the condition of coalgebra morphism for $i_{A}$)}\\=&i_{A}\circ\mu^{2}_{I(q_{\overline{A}})}\;\footnotesize\textnormal{(by \eqref{overB2mu})}.
\end{align*}
	
As a conclusion, due to the fact that $i_{A}$ is a monomorphism, $\overline{\mu}_{\IA}=\mu^{2}_{I(q_{\overline{A}})}$.
\end{proof}

\begin{theorem}
The diagram of functors
\[\xymatrix{
&{\sf P(coc\textnormal{-}HBr)}\ar[d]^-{{\sf P}'}\ar[rr]^-{{\sf Q}'} & &{\sf P(coc\textnormal{-}rRB)}\ar[d]^-{{\sf P}}\\
&{\sf coc\textnormal{-}HBr}\ar[rr]^-{{\sf F}} & &{\sf coc\textnormal{-}rRB}
}\]
is commutative.
\end{theorem}
\begin{proof}
Let $(\mathbb{H},\mathbb{D},x,y)\in{\sf P(coc\textnormal{-}HBr)}$. We have to show that the relative Rota-Baxter operator 
\begin{align*}
&({\sf P}\circ {\sf Q}')((\mathbb{H},\mathbb{D},x,y))\\=&{\sf P}\left(\left(\left(id_{H}\begin{array}{c}H_{1}\\\downarrow\\H_{2}\end{array},\Gamma_{H_{1}}\right),\left(id_{D}\begin{array}{c}D_{1}\\\downarrow\\D_{2}\end{array},\Gamma_{D_{1}}\right),x,x,y,y\right)\right)=\left(id_{\ID}\begin{array}{c}\ID\\\downarrow\\\ID\end{array},\varphi_{\ID}\right),
\end{align*}
where $\varphi_{\ID}$ is the unique action satisfying that 
\begin{equation}\label{PhiIBb0}
i_{D}\circ\varphi_{\ID}=\Gamma_{D_{1}}\circ (i_{D}\otimes i_{D}),
\end{equation}
is the same as the relative Rota-Baxter operator
\begin{align*}
&({\sf F}\circ{\sf P}')((\mathbb{H},\mathbb{D},x,y))={\sf F}(\mathbb{I}(q_{D}))=\left(id_{\ID}\begin{array}{c}\ID_{1}\\\downarrow\\\ID_{2}\end{array},\Gamma_{\ID_{1}}\right).
\end{align*}

To finish the proof it is enough to see that $\varphi_{\ID}=\Gamma_{\ID_{1}}$, what follows by \eqref{PhiIBb0}, by  the equality $i_{D}\circ \Gamma_{I(q_{D})_{1}}=\Gamma_{D_{1}}\circ (i_{D}\otimes i_{D})$ and by the fact that $i_{D}$ is a monomorphism.
\end{proof}

To finish the article we will introduce the notion of strong projection of relative Rota-Baxter operators in order to prove that any such projection gives rise to a module in the sense of Definition \ref{module_RB} in a cocommutative setting.

\begin{definition}
A projection of relative Rota-Baxter operators 
\[\left(\left(T\begin{array}{c}H\\\downarrow\\B \end{array},\varphi_{H}\right),\left(L\begin{array}{c}A\\\downarrow\\D \end{array},\varphi_{A}\right),f, h,g,l\right)\]
is said to be strong when the following conditions hold:
\begin{gather}
\label{condstrongmod}
p_{A}\circ\varphi_{A}=p_{A}\circ\varphi_{A}\circ (D\otimes q_{A}),\\
\label{strongPrRB}
\varphi_{\overline{A}}^{ad}\circ(f\otimes i_{A})=\mu_{A}\circ (\overline{\mu}_{A}\otimes\lambda_{A})\circ(A\otimes c_{A,A})\circ((\delta_{A}\circ f)\otimes i_{A}),
\end{gather} 
where $\overline{A}$ is the Hopf algebra introduced in Theorem \ref{adjunction}.

Strong projections of relative Rota-Baxter operators constitute a full subcategory of {\sf P(rRB)} which we will denote by ${\sf SP(rRB)}$. When the relative Rota-Baxter operators involved in the strong projection are cocommutative, they constitute a full subcategory of ${\sf SP(rRB)}$ denoted by ${\sf SP(coc\textnormal{-}rRB)}$.
\end{definition}

\begin{theorem}
If $$\left(\left(T\begin{array}{c}H\\\downarrow\\B \end{array},\varphi_{H}\right),\left(L\begin{array}{c}A\\\downarrow\\D \end{array},\varphi_{A}\right),f, h,g,l\right)$$ is an object in ${\sf SP(coc\text{-}rRB)}$, then 
$$(\IA,\ID,\psi_{\IA},\varphi_{\IA},\psi_{\ID},L_{0})$$
is a module over the relative Rota-Baxter operator $\left(T\begin{array}{c}H\\\downarrow\\B \end{array},\varphi_{H}\right)$, where $L_{0}\colon \IA\rightarrow \ID$ is the morphism defined by \eqref{T0cons} and the actions are defined as follows:
\begin{gather*}
\psi_{\IA}\coloneqq p_{A}\circ\mu_{A}\circ (f\otimes i_{A}),\\
\varphi_{\IA}\coloneqq p_{A}\circ \varphi_{A}\circ(h\otimes i_{A}),\\
\psi_{\ID}\coloneqq p_{D}\circ\mu_{C}\circ(h\otimes i_{D}).
\end{gather*} 
\end{theorem}
\begin{proof}
By the general theory of Hopf algebra projections, it is well-known that $(\IA,\psi_{\IA})$ is a left $H$-module and $(\ID,\psi_{\ID})$ is a left $B$-module. Let's show that $(\IA,\varphi_{\IA})$ is also a left $B$-module. Indeed, on the one hand it is straightforward to compute that $\varphi_{\IA}\circ(\eta_{B}\otimes\IA)=id_{\IA}$ and, on the other hand, we have that
\begin{align*}
&\varphi_{\IA}\circ(B\otimes\varphi_{\IA})\\=&p_{A}\circ\varphi_{A}\circ(h\otimes(q_{A}\circ\varphi_{A}\circ (h\otimes i_{A})))\;\footnotesize\textnormal{(by $q_{A}=i_{A}\circ p_{A}$)}\\=&p_{A}\circ\varphi_{A}\circ(h\otimes (\varphi_{A}\circ (h\otimes i_{A})))\;\footnotesize\textnormal{(by \eqref{condstrongmod})}\\=&p_{A}\circ\varphi_{A}\circ ((\mu_{D}\circ (h\otimes h))\otimes i_{A})\;\footnotesize\textnormal{(by module axioms for $(A,\varphi_{A})$)}\\=&\varphi_{\IA}\circ(\mu_{B}\otimes \IA)\;\footnotesize\textnormal{(by the condition of algebra morphism for $h$)}.
\end{align*}

To conclude the proof it only remains us to show that \eqref{compatmodHmodB} and \eqref{compatmodT} hold. At first, we have that
\begin{align*}
&\psi_{\IA}\circ(\varphi_{H}\otimes\varphi_{\IA})\circ (B\otimes c_{B,H}\otimes \IA)\circ(\delta_{B}\otimes H\otimes \IA)\\=&p_{A}\circ \mu_{A}\circ((f\circ\varphi_{H})\otimes (\varphi_{A}\circ (h\otimes i_{A})))\circ (B\otimes c_{B,H}\otimes\IA)\circ (\delta_{B}\otimes H\otimes \IA)\;\footnotesize\textnormal{(by \eqref{id1})}\\=&p_{A}\circ \mu_{A}\circ (\varphi_{A}\otimes\varphi_{A})\circ (D\otimes c_{D,A}\otimes A)\circ (((h\otimes h)\circ\delta_{B})\otimes f\otimes i_{A})\;\footnotesize\textnormal{(by \eqref{cond2morrRB} for $(f,h)$)}\\=&p_{A}\circ\mu_{A}\circ(\varphi_{A}\otimes\varphi_{A})\circ (D\otimes c_{D,A}\otimes A)\circ(\delta_{D}\otimes A\otimes A)\circ (h\otimes f\otimes i_{A})\;\footnotesize\textnormal{(by the condition of coalgebra}\\&\footnotesize\textnormal{morphism for $h$)}\\=&p_{A}\circ \varphi_{A}\circ (h\otimes(\mu_{A}\circ (f\otimes i_{A})))\;\footnotesize\textnormal{(by the condition of morphism of left $D$-modules for $\mu_{A}$)}\\=&\varphi_{\IA}\otimes(B\otimes \psi_{\IA})\;\footnotesize\textnormal{(by \eqref{condstrongmod})},
\end{align*}
what implies that \eqref{compatmodHmodB} holds. Let's define
\[\kappa\coloneqq L\circ \varphi_{\overline{A}}^{ad}\circ(f\otimes i_{A}).\]
On the one side, it is obtained that
\begin{align*}
&i_{D}\circ\psi_{\ID}\circ (T\otimes L_{0})\\=&\varphi_{D}^{ad}\circ((h\circ T)\otimes (i_{D}\circ L_{0}))\;\footnotesize\textnormal{(by \eqref{phi2})}\\=&\varphi_{D}^{ad}\circ ((L\circ f)\otimes(L\circ i_{A}))\;\footnotesize\textnormal{(by \eqref{T0cons} and \eqref{cond1morrRB} for $(f,h)$)}\\=&\kappa\;\footnotesize\textnormal{(by the condition of Hopf algebra morphism for $L\colon {\overline A}\rightarrow D$)},
\end{align*}
and, on the other side,
\begin{align*}
&i_{D}\circ L_{0}\circ \psi_{\IA}\circ (H\otimes(\varphi_{\IA}\circ (T\otimes \IA)))\circ(\delta_{H}\otimes\IA)\\=&L\circ q_{A}\circ\mu_{A}\circ(A\otimes\varphi_{A})\circ(((f\otimes (h\circ T))\circ \delta_{H})\otimes i_{A})\;\footnotesize\textnormal{(by \eqref{id1} and \eqref{T0cons})}\\=&L\circ\mu_{A}\circ (\mu_{A}\otimes(f\circ\lambda_{H}\circ g\circ\mu_{A}))\circ(A\otimes c_{A,A}\otimes A)\circ((\delta_{A}\circ f)\otimes (((\varphi_{A}\circ(L\otimes A))\otimes(\varphi_{A}\circ(L\otimes A)))\\&\circ (A\otimes c_{A,A}\otimes A)\circ ((\delta_{A}\circ f)\otimes(\delta_{A}\circ i_{A}))))\circ(\delta_{H}\otimes\IA)\;\footnotesize\textnormal{(by \eqref{cond1morrRB} for $(f,h)$ and the condition of coalgebra}\\&\footnotesize\textnormal{morphism for $\mu_{A}$, $\varphi_{A}$ and $L$)}\\=&L\circ\mu_{A}\circ (\mu_{A}\otimes(f\circ\lambda_{H}\circ \mu_{H}))\circ(A\otimes c_{H,A}\otimes H)\circ(((f\otimes H)\circ\delta_{H})\otimes (((\varphi_{A}\circ((L\circ f)\otimes A))\otimes(\varphi_{H}\\&\circ((l\circ L\circ f)\otimes H)))\circ (H\otimes c_{H,A}\otimes H)\circ (\delta_{H}\otimes((A\otimes g)\circ\delta_{A}\circ i_{A}))))\circ(\delta_{H}\otimes\IA)\;\footnotesize\textnormal{(by the condition of}\\&\footnotesize\textnormal{algebra morphism for $g$, \eqref{cond2morrRB} for $(g,l)$, naturality of $c$, the condition of coalgebra morphism for $f$ and $g\circ f=id_{H}$)}\\=&L\circ \mu_{A}\circ (\mu_{A}\otimes(f\circ\lambda_{H}))\circ(A\otimes c_{H,A})\circ(((f\otimes H)\circ\delta_{H})\otimes(\varphi_{A}\circ ((L\circ f)\otimes i_{A})))\circ(\delta_{H}\otimes\IA)\\&\footnotesize\textnormal{(by the equalizer condition for $i_{A}$, the fact that $\eta_{H}$ is a morphism of left $B$-modules and (co)unit properties)}\\=&L\circ\mu_{A}\circ(\overline{\mu}_{A}\otimes \lambda_{A})\circ(A\otimes c_{A,A})\circ(\delta_{A}\otimes A)\circ(f\otimes i_{A})\;\footnotesize\textnormal{(by \eqref{morant}, the condition of coalgebra morphism for $f$,}\\&\footnotesize\textnormal{naturality of $c$ and cocommutativity of $\delta_{A}$)}\\=&\kappa\;\footnotesize\textnormal{(by \eqref{strongPrRB})},
\end{align*}
what implies that \eqref{compatmodT} holds and the proof is concluded due to the fact that $i_{D}$ is a monomorphism.
\end{proof}


\section*{Funding Declaration}
The  authors were supported by  Ministerio de Ciencia e Innovaci\'on of Spain. Agencia Estatal de Investigaci\'on. Uni\'on Europea - Fondo Europeo de Desarrollo Regional (FEDER). Grant PID2020-115155GB-I00: Homolog\'{\i}a, homotop\'{\i}a e invariantes categ\'oricos en grupos y \'algebras no asociativas. Moreover, José Manuel Fernández Vilaboa and Brais Ramos Pérez were funded by Xunta de Galicia, grant ED431C 2023/31 (European FEDER support included, UE). Also, Brais Ramos Pérez was financially supported by Xunta de Galicia Scholarship ED481A-2023-023.



\bibliographystyle{amsalpha}

\begin{thebibliography}{A}

\bibitem{MN}  J. N. Alonso \'Alvarez and J. M. Fern\'andez Vilaboa, {\em Cleft extensions in braided categories}, \emph{Comm. Algebra} \textbf{28}(7) (2000)  3185-3196. 



\bibitem{AGV} I. Angiono, C. Galindo and L. Vendramin, Hopf braces and Yang-Baxter operators, \emph{Proc. Am. Math. Soc.} \textbf{145}(5) (2017) 1981-1995.

\bibitem{Baxter}  R. J. Baxter, Partition function of the eight-vertex lattice model, \emph{Ann. Phys.} \textbf{70}(1) (1972) 193-228.

\bibitem{Besp} Y. Bespalov, Crossed modules and quantum groups in braided categories, \emph{Appl. Categ. Structures} \textbf{5}(2) (1997) 155-204.

\bibitem{BCM} B. J. Blattner, M. Cohen, S. Montgomery, Crossed products and inner actions of Hopf algebras, \emph{Trans. Amer. Math. Soc.} \textbf{298}(2) (1986) 671-711.

\bibitem{BRZ1} T. Brzeziński, Trusses: Between braces and rings, \emph{Trans. Am. Math. Soc.} \textbf{372}(6) (2019) 4149-4176.

\bibitem{DR2} V. G. Drinfeld, Quantum groups, in \emph{Proc. Int. Cong. Math., (Berkeley, Calif. 1986)}, Vol. 1, 2, pp. 798-820.

\bibitem{DR1} V. G. Drinfeld, On some unsolved problems in quantum group theory, in \emph{Quantum groups (Leningrad, 1990)}, Vol. 1510 (Springer, Berlin, 1992), pp. 1-8.

\bibitem{ESS} P. Etingof, T. Schedler and A. Soloviev, Set-theoretical solutions to the quantum Yang-Baxter equation, \emph{Duke Math. J.} \textbf{100} (1999) 169-209.

\bibitem{FGR} J. M. Fernández Vilaboa, R. González Rodríguez and B. Ramos Pérez, Categorical isomorphisms for Hopf braces, Preprint (2023), arXiv:2311.05235.

\bibitem{FGRR} J. M. Fernández Vilaboa, R. González Rodríguez, B. Ramos Pérez and A. B. Rodríguez Raposo, Projections of Hopf braces, Preprint (2024), arXiv:2404.12231.

\bibitem{GI} T. Gateva-Ivanova, A combinatorial approach to the set-theoretic solutions of the Yang-Baxter equation, \emph{J. Math. Phys.} \textbf{45} (2004) 3828-3858.

\bibitem{Goncharov} M. Goncharov, Rota-Baxter operators on cocommutative Hopf algebras, \emph{J. Algebra} \textbf{582} (2021) 39-56.

\bibitem{RGON} R. González Rodríguez, The fundamental theorem of Hopf modules for Hopf braces, \emph{Linear Multilinear Algebra} \textbf{70}(20) (2022) 5146-5156.

\bibitem{GONROD} R. González Rodríguez and A. B. Rodríguez Raposo, Categorical equivalences for Hopf trusses and their modules, Preprint (2023), arXiv:2312.06520.

\bibitem{GV} L. Guarnieri and L. Vendramin, Skew braces and the Yang–Baxter equation, \emph{Math. Comput.} \textbf{86}(307) (2017) 2519-2534.




\bibitem{K} C. Kassel, {\em Quantum Groups} (Springer-Verlag, New York, 1995).

\bibitem{LST} Y. Li, Y. Sheng and R. Tang, Post-Hopf algebras, relative Rota-Baxter operators and solutions of the Yang-Baxter equation, \emph{J. Noncommut. Geom.} (2024) (in press: DOI 10.4171/JNCG/537).




\bibitem{MAJ2} S. Majid, Cross products by braided groups and bosonization, \emph{J. Algebra} \textbf{163}(1) (1994) 165-190. 

\bibitem{RAD} D. E. Radford, The structure of Hopf algebras with a projection, \emph{J. Algebra} \textbf{92}(2) (1985) 322-347.

\bibitem{Rump} W. Rump, Braces, radical rings, and the quantum Yang-Baxter equation, \emph{J Algebra} \textbf{307}(1) (2007) 153-170.

\bibitem{Sch} P. Schauenburg, On the braiding on a Hopf algebra in a braided category, \emph{New York J. Math} \textbf{4} (1998) 259-263.

\bibitem{SW} M. E. Sweedler,  {\em Hopf algebras} (Benjamin, New York 1969).


\bibitem{Yang} C. N. Yang, Some exact results for the many-body problem in one dimension with repulsive delta-function interaction, \emph{Phys. Rev. Lett.} \textbf{19}(23) (1967) 1312.

\end{thebibliography}

\end{document}